\documentclass[aos,seceqn,dvips]{arximspdf}
\usepackage{multirow}
\usepackage{dcolumn}
\usepackage{graphicx}
\doi{10.1214/10-AOS842}
\volume{39}
\issue{1}
\pubyear{2011}
\firstpage{305}
\lastpage{332}

\makeatletter
\newtheorem{thmm}{Theorem}[section]
\newtheorem{lemma}{Lemma}[section]
\newproclaim{example}{Example}
\newproclaim{remark}{Remark}

\newcolumntype{d}[1]{D{.}{.}{#1}}
\newcommand{\sumn}{\sum_{i=1}^{n}}

\newcommand{\sumnn}{\sum_{i'=1}^{n}}
\newcommand{\sumd}{\sum_{j=1}^{d}}
\newcommand{\sumdd}{\sum_{j=1}^{d_2}}

\newcommand{\sumq}{\sum_{k=1}^{q}}
\newcommand{\sumqq}{\sum_{k'=1}^{q}}

\newcommand{\nn}{\nonumber}

\newcommand{\diag}{\mathrm{diag}}
\newcommand{\argmin}{\operatorname{argmin}}

\newcommand{\Rb}{R_2(q)}
\newcommand{\wcon}{\stackrel{\mathcal{D}}{\longrightarrow}}
\newcommand{\pcon}{\stackrel{\mathcal{P}}{\longrightarrow}}
\newcommand{\LS}{\mathrm{LS}}
\newcommand{\ngrid}{n_{\mathrm{grid}}}
\newcommand{\rar}{\rightarrow}

\newcommand{\lam}{\lambda}
\newcommand{\alp}{\alpha}

\newcommand{\eps}{\varepsilon}

\newcommand{\bal}{{\bolds{\alpha}}}
\newcommand{\bbe}{{\bolds{\beta}}}
\newcommand{\bth}{{\bolds{\theta}}}
\newcommand{\bSig}{{\bolds{\Sigma}}}
\newcommand{\bDel}{{\bolds{\Delta}}}

\newcommand{\bXi}{{\bolds{\Xi}}}

\newcommand{\bA}{\mathbf{A}}
\newcommand{\bB}{\mathbf{B}}

\newcommand{\bS}{\mathbf{S}}

\newcommand{\bW}{\mathbf{W}}
\newcommand{\bX}{\mathbf{X}}
\newcommand{\bZ}{\mathbf{Z}}
\newcommand{\ba}{\mathbf{a}}
\newcommand{\bb}{\mathbf{b}}
\newcommand{\bc}{\mathbf{c}}
\newcommand{\be}{\mathbf{e}}

\newcommand{\bt}{\mathbf{t}}

\newcommand{\bx}{\mathbf{x}}

\newcommand{\bz}{\mathbf{z}}
\newcommand{\bzr}{\mathbf{0}}
\newcommand{\bone}{\mathbf{1}}
\newcommand{\bD}{\mathbf{D}}
\makeatother

\begin{document}
\begin{frontmatter}

\title{New efficient estimation and variable selection methods for semiparametric varying-coefficient
partially linear models}

\runtitle{Sparse semiparametric varying-coefficient PLMS}

\begin{aug}
\author[A]{\fnms{Bo} \snm{Kai}\thanksref{t1}\ead[label=e1]{KaiB@cofc.edu}},
\author[B]{\fnms{Runze} \snm{Li}\thanksref{t2}\ead[label=e2]{rli@stat.psu.edu}}
\and
\author[C]{\fnms{Hui} \snm{Zou}\thanksref{t3}\ead[label=e3]{hzou@stat.umn.edu}\corref{}}

\thankstext{t1}{Supported by NIDA, NIH Grants R21 DA024260 and P50 DA10075.}
\thankstext{t2}{Supported by NNSF of China Grant 11028103, NSF Grant DMS-03-48869 and  NIDA, NIH Grants R21 DA024260. The content is solely the responsibility of the authors and does not necessarily represent the official views of the NIDA or the NIH.}
\thankstext{t3}{Supported by NSF Grant DMS-08-46068.}
\runauthor{B. Kai, R. Li and H. Zou}
\affiliation{College of Charleston, Pennsylvania State University and University of Minnesota}
\address[A]{B. Kai\\
Department of Mathematics\\
College of Charleston\\
Charleston, South Carolina 29424\\
USA\\
\printead{e1}} %adresu isvedimo komanda gale!

\address[B]{R. Li\\
Department of Statistics\\
Pennsylvania State University\\
University Park, Pennsylvania 16802\\
USA\\
\printead{e2}}

\address[C]{H. Zou\\
School of Statistics\\
University of Minnesota\\
Minneapolis, Minnesota 55455\\
USA\\
\printead{e3}}
\end{aug}

% HISTORY:
\received{\smonth{1} \syear{2010}}
\revised{\smonth{5} \syear{2010}}

% ABSTRACT
\begin{abstract}
The complexity of semiparametric models poses new challenges to
statistical inference and model selection that frequently arise
from real applications. In this work, we propose new estimation
and variable selection procedures for the semiparametric
varying-coefficient partially linear model. We first study
quantile regression estimates for the nonparametric
varying-coefficient functions and the parametric regression
coefficients. To achieve nice efficiency properties, we further
develop a semiparametric composite quantile regression procedure.
We establish the asymptotic normality of proposed estimators for
both the parametric and nonparametric parts and show that the
estimators achieve the best convergence rate. Moreover, we show
that the proposed method is much more efficient than the
least-squares-based method for many non-normal errors and that it only
loses a small amount of efficiency for normal errors. In addition, it is
shown that the loss in efficiency is at most $11.1\%$ for
estimating varying coefficient functions and is no greater than
$13.6\%$ for estimating parametric components. To achieve sparsity
with high-dimensional covariates, we propose adaptive penalization
methods for variable selection in the semiparametric
varying-coefficient partially linear model and prove that the methods
possess the oracle property. Extensive Monte Carlo simulation
studies are conducted to examine the finite-sample performance of
the proposed procedures. Finally, we apply the new methods to
analyze the plasma beta-carotene level data.
\end{abstract}



% KEYWORDS

\end{frontmatter}

%s1 ###
\section{Introduction}\label{s1}

Semiparametric regression modeling has recently become popular in the
statistics literature because it keeps the flexibility of nonparametric
models while maintaining the explanatory power of parametric models.
The partially linear model, the most commonly used semiparametric
regression model, has received a lot of attention in the literature;
see H{\"a}rdle, Liang and Gao \cite{HardleLiangGao2000}, Yatchew
\cite{Yatchew2003} and references therein for theory and applications
of partially linear models. Various extensions of the partially linear
model have been proposed in the literature; see Ruppert, Wand and
Carroll \cite{RuppertWandCarroll2003} for applications and theoretical
developments of semiparametric regression models. The semiparametric
varying-coefficient partially linear model, as an important extension
of the partially linear model, is becoming popular in the literature.
Let $Y$ be a response variable and $\{U, \bX, \bZ\}$ its covariates.
The semiparametric varying-coefficient partially linear model is
defined to be
%e1.1 ###
\begin{equation}\label{semiplm}
Y=\alp_0(U)+\bX^T\bal(U)+\bZ^T\bbe+\eps,
\end{equation}
where $\alp_0(U)$ is a baseline\vspace*{1.5pt} function,
$\bal(U)=\{\alp_1(U),\ldots,\alp_{d_1}(U)\}^T$ consists of $d_1$
unknown varying coefficient functions,
$\bbe=(\beta_1,\ldots,\beta_{d_2})^T$ is a $d_2$-dimensional
coefficient vector and $\eps$ is random error. In this paper, we will
focus on univariate $U$ only, although the proposed procedure is
directly applicable for multivariate~$\mathbf{U}$. Zhang, Lee and Song
\cite{ZhangLeeSong2002} proposed an estimation procedure for the model
(\ref{semiplm}), based on local polynomial regression techniques. Xia,
Zhang and Tong \cite{XiaZhangTong2004} proposed a semilocal estimation
procedure to further reduce the bias of the estimator for $\bbe$
suggested in Zhang, Lee and Song \cite{ZhangLeeSong2002}. Fan and Huang
\cite{FanHuang2005} proposed a profile least-squares estimator for
model (\ref{semiplm}) and developed statistical inference procedures.
As an extension of Fan and Huang \cite{FanHuang2005}, a profile
likelihood estimation procedure was developed in Lam and Fan
\cite{LamFan2008}, under the generalized linear model framework with a
diverging number of covariates.

Existing estimation procedures for model (\ref{semiplm}) were built on
either least-squares- or likelihood-based methods. Thus, the existing
procedures are expected to be sensitive to outliers and their
efficiency may be significantly improved for many commonly used
non-normal errors. In this paper, we propose new estimation procedures
for model (\ref{semiplm}). This paper contains three major
developments: (a)~semiparametric quantile regression; (b)~semiparametric composite quantile regression; (c)~adaptive penalization
methods for achieving sparsity in semiparametric composite quantile
regression.

Quantile regression is often considered as an alternative to
least-squares in the literature. For a complete review on quantile
regression, see Koenker \cite{Koenker2005}. Quantile-regression-based
inference procedures have been considered in the literature; see, for
example, Cai and Xu \cite{CaiXu08}, He and Shi \cite{HeShi1996}, He,
Zhu and Fung \cite{HeZhuFung2002}, Lee \cite{Lee2003}, among others. In
Section \ref{s2}, we propose a new semiparametric quantile regression
procedure for model (\ref{semiplm}). We investigate the sampling
properties of the proposed method and their asymptotic normality. When
applying semiparametric quantile regression to model (\ref{semiplm}),
we observe that all quantile regression estimators can estimate
$\bal(u)$ and $\bbe$ with the optimal rate of convergence. This fact
motivates us to combine the information across multiple quantile
estimates to obtain improved estimates of $\bal(u)$ and $\bbe$. Such an
idea has been studied for the parametric regression model in Zou and
Yuan \cite{ZouYuan2008} and it leads to the composite quantile
regression (CQR) estimator that is shown to enjoy nice asymptotic
efficiency properties compared with the classical least-squares
estimator. In Section~\ref{s3}, we propose the semiparametric composite
quantile regression (semi-CQR) estimators for estimating both
nonparametric and parametric parts in model (\ref{semiplm}). We show
that the semi-CQR estimators achieve the best convergence rates. We
also prove the asymptotic normality of the semi-CQR estimators. The
asymptotic theory shows that, compared with the semiparametric
least-squares estimators, the semi-CQR estimators can have substantial
efficiency gain for many non-normal errors and only lose a small amount
of efficiency for normal errors. Moreover, the relative efficiency is
at least $88.9\%$ for estimating varying-coefficient functions and is
at least $86.4\%$ for estimating parametric components.

In practice, there are often many covariates in the parametric part of
model~(\ref{semiplm}). With high-dimensional covariates, sparse
modeling is often considered superior, owing to enhanced model
predictability and interpretability \cite{FanLi2006}. Variable
selection for model~(\ref{semiplm}) is challenging because it involves
both nonparametric and parametric parts. Traditional variable selection
methods, such as stepwise regression or best subset variable selection,
do not work effectively for the semiparametric model because they need
to choose smoothing parameters for each submodel and cannot cope with
high-dimensionality. In Section \ref{s4}, we develop an effective
variable selection procedure to select significant parametric
components in model~(\ref{semiplm}). We demonstrate that the proposed
procedure possesses the oracle property, in the sense of Fan and Li
\cite{FanLi2001}.

In Section \ref{s5}, we conduct simulation studies to examine the
finite-sample performance of the proposed procedures. The proposed
methods are illustrated with the plasma beta-carotene level data.
Regularity conditions and technical proofs are given in Section
\ref{s6}.

%s2 ###
\section{Semiparametric quantile regression}\label{s2}

In this section, we develop the semiparametric quantile regression
method and theory. Let $\rho_{\tau}(r) = \tau r - r I$ $(r<0)$ be the
check loss function at $\tau\in(0,1)$. Quantile regression is often
used to estimate the conditional quantile functions of $Y$,
\[
Q_{\tau}(u,\bx,\bz)=\argmin\limits_{a}E\{\rho_{\tau}(Y-a)|(U,\bX,\bZ)=(u,\bx,\bz)\}.
\]
The semiparametric varying-coefficient partially linear model assumes
that the conditional quantile function is expressed as
$Q_{\tau}(u,\bx,\bz)=\alp_{0,\tau}(u)+\bx^T\bal_{\tau}(u)+\bz^T\bbe_{\tau}.$

Suppose that $\{U_i,\bX_i,\bZ_i,Y_i\}$, $i=1,\ldots,n,$ is an
independent and identically distributed sample from the model
%e2.1 ###
\begin{equation}\label{semiplmQr}
Y=\alp_{0,\tau}(U)+\bX^T\bal_{\tau}(U)+\bZ^T\bbe_{\tau}+\eps_{\tau},
\end{equation}
where $\eps_{\tau}$ is random error with conditional $\tau$th quantile
being zero. {We obtain quantile regression estimates of $\alp_{0,
\tau}(\cdot)$, $\bal_{\tau}(\cdot)$ and $\bbe_{\tau}$ by minimizing the
quantile loss function}
%e2.2 ###
\begin{equation}\label{quantile}
\sumn \rho_{\tau}\{Y_i-\alp_{0}(U_i)-\bX_i^T\bal(U_i)-\bZ_i^T\bbe\}.
\end{equation}
{Because (\ref{quantile}) involves both nonparametric and parametric
components, and because they can be estimated by different rates of
convergence, we propose a three-stage estimation procedure}. In the
first stage, we employ local linear regression techniques to derive an
initial estimates of $\alpha_{0,\tau}(\cdot)$, $\bal_{\tau}(\cdot)$
and $\bbe_{\tau}$. Then, in the second and third stages, we further
improve the estimation efficiency of the initial estimates for
$\bbe_{\tau}$ and $(\alpha_{0,\tau}(\cdot),\bal_{\tau}(\cdot))$,
respectively.

For $U$ in the neighborhood of $u$, we use a local linear approximation
\[
\alp_j(U)\approx\alp_j(u)+\alp'_j(u)(U-u)\triangleq a_j+b_j(U-u)
\]
for $j=0,\ldots,d_1$. Let $\{{\tilde{a}}_{0,\tau},{\tilde{b}}_{0,\tau},{\tilde\ba}_{\tau},{\tilde{\bb}}_\tau,{\tilde{\bbe}}_{\tau}\}$ be the minimizer
of the local weighted quantile loss function
\[%\label{elltau}
%&&\{\tilde{a}_{0, \tau}, \tilde{b}_{0, \tau}, \tilde\ba_{\tau}, \tilde\bb_\tau, \tilde\bbe_{\tau}\}\nn\\
%%    &&\ell_{\tau}(a_0,b_0,\ba,\bb,\bbe) \nn \\
%    &= & \argmin_{\{a_0,b_0,\ba,\bb,\bbe\}}
\sumn\rho_{\tau}\bigl\{Y_i-a_0-b_0(U_i-u)-\bX_i^T\{\ba+\bb(U_i-u)\}-\bZ_i^T\bbe\bigr\}K_h(U_i-u),
\]
where $\ba=(a_1,\ldots,a_{d_1})^T$, $\bb=(b_1,\ldots,b_{d_1})^T$,
$K(\cdot)$ is a given kernel function and $K_h(\cdot)=K(\cdot /h) / h$
{with a bandwidth $h$}. Then,
\[
{\tilde{\alp}}_{0,\tau}(u)={\tilde{a}}_0,\qquad{\tilde{\bal}}_{\tau}(u)={\tilde{\ba}}_{\tau}.
\]
%Denote the minimizer of (\ref{elltau}) by $\{\tilde{a}_{0, \tau}, \tilde{b}_{0, \tau},
{We take $\{\tilde{\alp}_{0,\tau}(u),\tilde\bal_{\tau}(u),\tilde\bbe_{\tau}\}$ as the initial estimates.

We now provide theoretical justifications for the initial estimates.
First, we give some notation. Let $f_\tau(\cdot|u,\bx,\bz)$ and
$F_\tau(\cdot|u,\bx,\bz)$  be the density function and cumulative
distribution function of the error conditional on
$(U,\bX,\bZ)=(u,\bx,\bz)$, respectively. Denote by $f_U(\cdot)$ the
marginal density function of the covariate $U$. The kernel $K(\cdot)$
is chosen as a symmetric density function and we let
\[
\mu_j=\int u^j K(u)\,du\quad\mbox{and}\quad\nu_j=\int u^j K^2(u)\,du,\qquad j=0,1,2,\ldots.
\]
%Let $\bbe_0$ be the true value of $\bbe$.
We then have the following result.

\begin{thmm}\label{ThmQrAlpha1}
Under the regularity conditions given in Section \textup{\ref{s7}}, if $h\rar0$
and $nh\rar\infty$ as $n\rar\infty$, then
\begin{eqnarray}\label{asymnorm41}
&&
\sqrt{nh}\left[
\pmatrix{
\tilde\alp_{0,\tau}(u)-\alp_{0,\tau}(u)\vspace*{3pt}\cr
\tilde \bal_{\tau}(u) - \bal_{\tau}(u)\vspace*{3pt}\cr
\tilde\bbe_{\tau} - \bbe_{\tau}
}-
\frac{\mu_2 h^2}{2}
\pmatrix{
 \alp''_{0, \tau}(u)\vspace*{3pt}\cr
 \bal''_{\tau}(u)\vspace*{2pt}\cr
 \bzr}\right]\nn
%$\tilde \bth^* = \sqrt{nh}(\tilde a_{01}-\alp_0(u)-c_1, \cdots, \tilde a_{0q}-\alp_0(u)-c_q, (\tilde \ba-\bal(u))^T, (\tilde \bbe-\bbe_0)^T, h (\tilde  b_0 - \alp'_0(u) ), h (\tilde  \bb - \bal'(u) )^T)^T$
\\[-8pt]\\[-8pt]
&&\qquad\wcon N\biggl(\bzr,\frac{\nu_0\tau(1-\tau)}{f_U(u)%f^2(c_{\tau})
}\bA_1^{-1}(u)\bB_1(u)\bA_1^{-1}(u)\biggr),\nn
\end{eqnarray}
where
%$c_{\tau} = F^{-1}(\tau)$, $\alp_{0, \tau}(u)=\alp_0(u) + c_{\tau}$,
$\bA_1(u)=E[
f_\tau(0|U,\bX,\bZ)(1,\bX^T,\bZ^T)^T(1,\bX^T,\bZ^T)|U=u ]$
 and
$\bB_1(u)=E[
(1,\bX^T,\bZ^T)^T(1,\bX^T,\bZ^T)|U=u ]$.
\end{thmm}

Theorem~\ref{ThmQrAlpha1} implies that $\tilde\bbe_\tau$ is a
$\sqrt{nh}$-consistent estimator---this is because we only use data in
a local neighborhood of $u$ to estimate\vspace*{-2pt} $\bbe_{\tau}$. Define
$Y_{i,\tau}^*=Y_i-\tilde\alp_{0,\tau}(U_i)-\bX_i^T\tilde\bal_\tau(U_i)$ and compute an improved estimator of $\bbe_{\tau}$ by
%e2.3 ###
\begin{equation}\label{elltau2}
\hat{\bbe}_{\tau}
=
\argmin\limits_{\bbe}\sumn\rho_{\tau}(Y_{i,\tau}^*-\bZ_i^T\bbe).
\end{equation}
We call it the \textit{semi-QR estimator} of $\bbe_{\tau}$. The next
theorem shows the asymptotic properties of  $\hat \bbe_{\tau}$.

\begin{thmm}\label{ThmQrBeta1} Let
$\bolds{\xi}_{\tau}(u,\bx,\bz)=E[f_{\tau}(0|U,\bX,\bZ)\bZ(1,\bX^T,\mathbf{0})|U=u]\times\break\bA_1^{-1}(u)(1,\bx^T,\bz^T)^T$. Under the
regularity conditions given in Section \textup{\ref{s7}}, if $nh^4\rightarrow 0$ and $nh^2/\log(1/h)\rightarrow\infty$ as $n\rightarrow\infty$,
then the asymptotic distribution of $\hat\bbe_{\tau}$ is given by
%e2.4 ###
\begin{equation}
\sqrt{n}(\hat\bbe_{\tau}-\bbe_{\tau})\wcon N(0,
    %\tau(1-\tau)
%    \frac{\tau(1-\tau)}{f^2(c_{\tau})}
\bS_{\tau}^{-1}\bXi_{\tau}\bS_{\tau}^{-1}),
\end{equation}
where $\bS_{\tau}=E[f_\tau(0|U,\bX,\bZ)\bZ\bZ^T]$ and
$\bXi_{\tau}=\tau(1-\tau)E[\{\bZ-\bolds{\xi}_{\tau}(U,\bX,\bZ)\}\{\bZ-\bolds{\xi}_{\tau}(U,\bX,\bZ)\}^T]$.
\end{thmm}

The optimal bandwidth in Theorem~\ref{ThmQrAlpha1} is $h\sim n^{-1/5}$. This bandwidth does not satisfy the condition in
Theorem~\ref{ThmQrBeta1}. Hence, in order to obtain the root-$n$
consistency and asymptotic normality for $\hat\bbe_{\tau}$,
undersmoothing for $\tilde{\alpha}_{0,\tau}(u)$ and
$\tilde{\bal}_{\tau}(u)$ is necessary. This is a common requirement in
semiparametric models; see Carroll et~al. \cite{CarrollEtal1997} for a
detailed discussion.

After obtaining the root-$n$ consistent estimator $\hat\bbe_{\tau}$,
we can further improve the efficiency of $\tilde\alp_{0,\tau}(u)$ and
$\tilde\bal_{\tau}(u)$. To this end, let $\{\hat{a}_{0,\tau},\hat{b}_{0,\tau},\hat{\ba}_{\tau},\hat{\bb}_{\tau}\}$ be the
minimizer of
%$Y_{i,\tau}^{**}= Y_i  - \bZ_i^T\hat\bbe_{\tau}$, and
\[
%=\argmin_{\{a_{0}, b_{0}, \ba,\bb\}}
\sumn\rho_{\tau}\bigl\{Y_i-\bZ_i^T \hat\bbe_{\tau}-a_{0}-b_0(U_i-u)-\bX_i^T \{\ba+\bb(U_i-u)\}\bigr\}K_h(U_i-u).
\]
%Denote the minimizer of (\ref{ell3}) by $\{\hat \ba_0,\hat b_0,\hat \ba,\hat \bb\}$, then
We define
%e2.5 ###
\begin{equation}\label{RefinedAplha42}
\hat\alp_{0,\tau}(u)=\hat{a}_{0,\tau},
\qquad
\hat\bal_\tau(u)=\hat\ba_{\tau}.
\end{equation}

\begin{thmm} \label{ThmQrAlpha2}
Under the regularity conditions given in Section \textup{\ref{s7}}, if $h\rar0$
and $nh\rar\infty$ as $n\rar\infty$, then
\begin{eqnarray}\label{asymnorm3}
&&
\sqrt{nh}\left[\pmatrix{
\hat\alp_{0,\tau}(u)-\alp_{0,\tau}(u)\vspace*{3pt}\cr
\hat\bal_{\tau}(u)-\bal_{\tau}(u)
}-
\frac{\mu_2 h^2}{2}
\pmatrix{
\alp''_{0,\tau}(u)\vspace*{3pt}\cr
\bal''_{\tau}(u)
}
%$\tilde \bth^* = \sqrt{nh}(\tilde a_{01}-\alp_0(u)-c_1, \cdots, \tilde a_{0q}-\alp_0(u)-c_q, (\tilde \ba-\bal(u))^T, (\tilde \bbe-\bbe_0)^T, h (\tilde  b_0 - \alp'_0(u) ), h (\tilde  \bb - \bal'(u) )^T)^T$
\right]\nn
\\[-8pt]\\[-8pt]
&&\qquad\wcon
N\biggl(\bzr,\frac{\nu_0\tau(1-\tau)}{f_U(u)
%f^2(c_{\tau}) } \bB^{-1}(u)
}\bA_2^{-1}(u)\bB_2(u)\bA_2^{-1}(u)\biggr),\nn
\end{eqnarray}
where $\bA_2(u)=E[f_{\tau}(0|U,\bX,\bZ)(1,\bX^T)^T(1,\bX^T)|U=u]$ and $\bB_2(u)=E[(1,\break\bX^T)^T(1,\bX^T)|U=u ]$.
\end{thmm}

Theorem~\ref{ThmQrAlpha2} shows that $\hat\alp_{0, \tau}(u)$ and $\hat
\bal_{\tau}(u)$ have the same conditional asymptotic biases as
$\tilde\alp_{0, \tau}(u)$ and $\tilde \bal_{\tau}(u)$, while they have
smaller conditional asymptotic variances than $\tilde\alp_{0, \tau}(u)$
and $\tilde \bal_{\tau}(u)$, respectively. Hence, they are
asymptotically more efficient than $\tilde{\alpha}_{0,\tau}(u)$ and
$\tilde{\bal}_{\tau}(u)$.

%s3 ###
\section{Semiparametric composite quantile regression}\label{s3}

The analysis of semiparametric quantile regression in Section \ref{s2}
provides a solid foundation for developing the semiparametric composite
quantile regression (CQR) estimates. We consider the connection between
the quantile regression model~(\ref{semiplmQr}) and model~(\ref{semiplm}) in the situations where the random error $\eps$ is
independent of $(U,\bX,\bZ)$. Let us assume that $Y=\alp_{0}(U)+\bX^T\bal(U)+\bZ^T\bbe+\eps$, where $\eps$ follows a distribution
$F$ with mean zero. In such situations, $Q_\tau(u,\bx,\bz)=\alpha_0(u)+c_\tau+\bx^T\bal(u)+\bz^T\bbe$, where
$c_{\tau}=F^{-1}(\tau)$. Thus, all quantile regression estimates
[$\hat{\bal}_{\tau}(u)$ and $\hat{\bbe}_{\tau}$ for all $\tau$]
estimate the same target quantities [$\bal(u)$ and $\bbe$] with the
optimal rate of convergence. Therefore, we can consider combining the
information across multiple quantile estimates to obtain improved
estimates of $\bal(u)$ and $\bbe$. Such an idea has been studied for
the parametric regression model, in Zou and Yuan \cite{ZouYuan2008},
and it leads to the CQR estimator that is shown to enjoy nice
asymptotic efficiency properties compared with the classical
least-squares estimator. Kai, Li and Zou \cite{KaiLiZou2009a} proposed
the local polynomial CQR estimator for estimating the nonparametric
regression function and its derivative. It is shown that the local CQR
method can significantly improve the estimation efficiency of the local
least-squares estimator for commonly used non-normal error
distributions. Inspired by these nice results, we study semiparametric
CQR estimates for model (\ref{semiplm}).

Suppose $\{U_i,\bX_i,\bZ_i,Y_i,i=1,\ldots,n\}$ is an independent and
identically distributed sample from model (\ref{semiplm}) and $\eps$
has mean zero. For a given $q$, let $\tau_k=k/(q+1)$ for $
k=1,2,\ldots,q$. The CQR procedure estimates $\alp_0(\cdot)$,
$\bal(\cdot)$ and $\bbe$ by minimizing the CQR loss function,
\[%\label{semicqr}
\sumq\sumn\rho_{\tau_k}\{Y_i-\alp_{0k}(U_i)-\bx_i^T\bal(U_i)-\bz_i^T\bbe\}.
\]
To this end, we adapt the three-stage estimation procedure from Section
\ref{s2}.

First, we derive good initial semi-CQR estimates. Let $\{\tilde\ba_0,
\tilde b_0,\tilde\ba,\tilde\bb,\tilde\bbe\}$ be the minimizer of the
local CQR loss function
\[%\label{ell1cqr}
\sumq\sumn\rho_{\tau_k}\bigl\{Y_i-a_{0k}-b_0(U_i-u)-\bX_i^T \{\ba+\bb(U_i-u)\}-\bZ_i^T\bbe\bigr\}K_h(U_i-u),
\]
where $\ba_0=(a_{01},\ldots,a_{0q})^T$, $\ba=(a_1,\ldots,a_{d_1})^T$
and $\bb=(b_1,\ldots,b_{d_1})^T$. Initial estimates of $\alp_{0}(u)$
and $\bal(u)$ are then given by
%e3.1 ###
\begin{equation}
\tilde\alp_{0}(u)=\frac{1}{q}\sumq\tilde{a}_{0k},\qquad\tilde\bal(u)=\tilde\ba.
\end{equation}

To investigate asymptotic behaviors of $\tilde\alp_0(u)$, $\tilde\bal(u)$ and $\tilde\bbe$, let us begin with some new {notation.}
Denote by $f(\cdot)$ and $F(\cdot)$  the density function and
cumulative distribution function of the error, respectively. Let $c_k=F^{-1}(\tau_k)$ and $C$ be a $q\times q$ diagonal matrix with
$C_{jj}=f(c_j)$. Write $\bc=C\bone$, $c=\bone^T C\bone$ and
\[
\bD_1(u)
=
E\left[\pmatrix{
C &\bc\bX^T&\bc\bZ^T\vspace*{2pt}\cr
\bX\bc^T&c\bX\bX^T&c\bX\bZ^T\vspace*{2pt}\cr
\bZ\bc^T&c\bZ\bX^T&c\bZ\bZ^T
}\biggl| U=u\right].
\]
Let $\tau_{kk'}=\tau_{k}\wedge\tau_{k'}-\tau_{k}\tau_{k'}$ and let
$T$ be a $q\times q$ matrix with the  $(k,k')$ element being
$\tau_{kk'}$. Write $\bt=T\bone$, $t=\bone^T T\bone$ and
\[
\bSig_1(u)
= E
\left[\pmatrix{
T&\bt\bX^T&\bt\bZ^T\vspace*{2pt}\cr
\bX\bt^T&t\bX\bX^T&t\bX\bZ^T\vspace*{2pt}\cr
\bZ\bt^T&t\bZ\bX^T&t\bZ\bZ^T
}\biggl|U=u\right].
\]

The following theorem describes the asymptotic sampling distribution of
$\{\tilde\ba_0,\tilde b_0,\tilde\ba,\tilde\bb,\tilde\bbe\}$.

\begin{thmm}\label{ThmCqrAlpha1}
Under the regularity conditions given in Section \ref{s7}, if $h\rar 0$ and $nh\rar\infty$ as $n\rar\infty$, then
\begin{eqnarray*}%\label{asymnorm45}
&&
\sqrt{nh}\left[
\pmatrix{
\tilde\ba_{0}-\bal_{0}(u)\vspace*{2pt}\cr
\tilde\ba-\bal(u)\vspace*{2pt}\cr
\tilde\bbe-\bbe_0
}
-\frac{\mu_2 h^2}{2}
\pmatrix{
\bal''_0(u)\vspace*{3pt}\cr
\bal''(u)\vspace*{2pt}\cr
\bzr
}
%$\tilde \bth^* = \sqrt{nh}(\tilde a_{01}-\alp_0(u)-c_1, \cdots, \tilde a_{0q}-\alp_0(u)-c_q, (\tilde \ba-\bal(u))^T, (\tilde \bbe-\bbe_0)^T, h (\tilde  b_0 - \alp'_0(u) ), h (\tilde  \bb - \bal'(u) )^T)^T$
\right]
\\
&&\qquad\wcon
N\biggl(\bzr,\frac{\nu_0}{f_U(u)}\bD^{-1}_1(u)\bSig_1(u)\bD^{-1}_1(u)\biggr),
\end{eqnarray*}
where $\bal_{0}(u)=(\alp_0(u)+c_1,\ldots,\alp_0(u)+c_q)^T$ and
$\bbe_0$ is the true value of $\bbe$.
\end{thmm}

With the initial estimates in hand, we are now ready to derive a
$\sqrt{n}$-consistent estimator of $\bbe$ by
%e3.2 ###
\begin{equation}\label{betacheck2}
\hat\bbe=\argmin\limits_{\bbe}\sumq\sumn\rho_{\tau_k}\{Y_i-\tilde a_{0k}(U_i)-\bX_i^T\tilde\ba(U_i)-\bZ_i^T\bbe\},
\end{equation}
which is called the \textit{semi-CQR estimator} of $\bbe$.

\begin{thmm}\label{ThmCqrBeta1}
Under the regularity conditions given in Section \ref{s7}, if $nh^4\rightarrow 0$ and $nh^2/\log(1/h)\rightarrow\infty$ as
$n\rightarrow\infty$, then the asymptotic distribution of $\hat\bbe$ is given by
%e3.3 ###
\begin{equation}
\sqrt{n}(\hat\bbe-\bbe_0)
\wcon
N\biggl(0,\frac{1}{c^2}\bS^{-1}\bDel\bS^{-1}\biggr),
\end{equation}
where $\bS=E(\bZ\bZ^T)$ and $\bDel=\sum_{k=1}^q\sum_{k'=1}^q\tau_{kk'}E[\{\bZ-\bolds{\delta}_{k}
(U,\bX,\bZ)\}\{\bZ-\bolds{\delta}_{k'}(U,\bX,\break\bZ)\}^T]$, with
$\bolds{\delta}_{k}(u,\bx,\bz)$ being the $k$th column of the
$d_2\times q$ matrix
\[
\bolds{\delta}(u,\bx,\bz)
=
E[\bZ(\bc^T,c\bX^T,\mathbf{0})|U=u]\bD_1^{-1}(u)(I_q,\bone^T\bx,\bone^T\bz)^T.
\]
\end{thmm}

Finally, $\hat\bbe$ can also be used to further refine the
estimates for the nonparametric part.} %For $\alp_0(u)$ and $\bal(u)$,
Let $\{\hat\ba_0,\hat b_0,\hat\ba,\hat\bb\}$  be the minimizer of
\[%\label{ell3}
%&&\{\hat \ba_0,\hat b_0,\hat \ba,\hat \bb\}  = \argmin_{\ba_0,b_0,\ba,\bb}\\
%&=&\argmin_{\ba_0,b_0,\ba,\bb}
\sumq\sumn\rho_{\tau_k}[Y_i-\bZ_i^T\hat\bbe-a_{0k}-b_0(U_i-u)-\bX_i^T \{\ba+\bb(U_i-u)\}]K_h(U_i-u),
\]
where $\ba_0=(a_{01},\ldots,a_{0q})^T$. We then define the
\textit{semi-CQR estimators} for $\alp_0(u)$ and $\bal(u)$ as
%e3.4 ###
\begin{equation}\label{RefinedAplha47}
\hat\alp_0(u)=\frac{1}{q}\sumq \hat a_{0k},\qquad\hat\bal(u)=\hat\ba.
\end{equation}

We now study the asymptotic properties of $\hat\alp_0(u)$ and $\hat\bal(u)$. Let
\begin{eqnarray*}
\bD_2(u)
&=&
E\left[\pmatrix{
C&\bc\bX^T\vspace*{2pt}\cr
\bX\bc^T&c\bX\bX^T
}\biggl|U=u\right],
\\
\bSig_2(u)
&=&
E\left[\pmatrix{
T&\bt\bX^T\vspace*{2pt}\cr
\bX \bt^T & t\bX\bX^T
}\biggl|U=u\right].
\end{eqnarray*}

\begin{thmm}\label{ThmCqrAlpha}
Under the regularity conditions given in Section \textup{\ref{s7}}, if $h\rar 0$
and $nh\rar\infty$ as $n\rar\infty$,  the asymptotic distributions of
$\hat\alp_0(u)$ and $\hat\bal(u)$ are given by
\begin{eqnarray*}
&&
\sqrt{nh}\Biggl(\hat\alp_{0}(u)-\alp_{0}(u)-\frac{1}{q}\sumq c_k-\frac{\mu_2 h^2}{2}\alp''_0(u)\Biggr)
\\
&&\qquad\wcon
N\biggl(0,\frac{\nu_0}{f_U(u)}\frac{1}{q^2}\bone^T[\bD_2^{-1}(u)\bSig_2(u)\bD_2^{-1}(u)]_{11}\bone\biggr)
\end{eqnarray*}
and
\[
\sqrt{nh}\biggl(\hat\bal(u)-\bal(u)-\frac{\mu_2 h^2}{2}\bal''(u)\biggr)
\wcon
N\biggl(0,\frac{\nu_0}{f_U(u)}[\bD_2^{-1}(u)\bSig_2(u)\bD_2^{-1}(u)]_{22}\biggr),
\]
where $[\cdot]_{11}$ denotes the upper-left $q\times q$ submatrix and
$[\cdot]_{22}$ denotes the lower-right $d_1\times d_1$ submatrix.
\end{thmm}

\begin{remark}\label{rem1}
$\bal(u)$ and $\bbe$ represent the contributions of covariates. They
are the central quantities of interest in semiparametric inference. Li
and Liang \cite{LiLiang2008} studied the least-squares-based
semiparametric estimation, which we will refer to as ``semi-LS'' in
this work. The major advantage of semi-CQR over the classical semi-LS
is that semi-CQR has competitive asymptotic efficiency. Furthermore,
semi-CQR is also more stable and robust. Intuitively speaking, these
advantages come from the fact that semi-CQR utilizes information shared
across multiple quantile functions, whereas semi-LS only uses the
information contained in the mean function.
\end{remark}

To elaborate on Remark~\ref{rem1}, we discuss the relative efficiency
of semi-CQR relative to semi-LS. Note that $E(Y|U)=\alpha_0(U)+E(\bX|U)^T\bal(U)+E(\bZ|U)^T\bbe.$ It then follows that
$Y=E(Y|U)+\{\bX-E(\bX|U)\}^T\bal(U)+\{\bZ-E(\bZ|U)\}^T\bbe+\epsilon$.
Without loss of generality, let us consider the situation in which
$E(\bX|U)=0$ and $E(\bZ|U)=0$. Then, all $\bD_1(u),\bD_2(u),\bSig_1(u)$ and $\bSig_2(u)$ become block diagonal matrices. Thus, from
Theorem~\ref{ThmCqrAlpha}, we have
\[
\sqrt{nh}\Biggl(\hat\alp_{0}(u)-\alp_{0}(u)-\frac{1}{q}\sumq c_k-\frac{\mu_2 h^2}{2}\alp''_0(u)\Biggr)
\wcon
N\biggl(0,R_1(q)\frac{\nu_0}{f_U(u)}\biggr)
\]
and
\[
\sqrt{nh}\biggl(\hat\bal(u)-\bal(u)-\frac{\mu_2 h^2}{2}\bal''(u)\biggr)
\wcon
N\biggl(0,R_2(q)\frac{\nu_0}{f_U(u)}E^{-1}(\bX\bX^T|U=u)\biggr),
\]
where
\[
R_1(q) = \frac{1}{q^2} \bone^T C^{-1} T C^{-1} \bone =
\frac{1}{q^2}\sum_{k=1}^q\sum_{k'=1}^q
\frac{\tau_{kk'}}{f(c_k)f(c_{k'})}
\]
and
\[
R_2(q)=\frac{t}{c^2}=\frac{\sum_{k=1}^q\sum_{k'=1}^q\tau_{kk'}}{\{\sum_{k=1}^q f(c_k)\}^2}.
\]
Note that
\[
\bolds{\delta}(u,\bx,\bz)=E(\bZ\bX^T,\bzr|U=u)
E\left[\pmatrix{
\bX\bX^T&\bX\bZ^T\vspace*{2pt}\cr
\bZ\bX^T&\bZ\bZ^T
}\biggl|U=u\right]^{-1}(\bone^T\bx,\bone^T\bz)^T
\]
with all columns of $\bolds{\delta}(u,\bx,\bz)$ the same. Thus, $\bDel=t\bDel_0$ with
$\bDel_0=E[\{\bZ-\bolds{\delta}_1(U,\bX,\bZ)\}\{\bZ-\bolds{\delta}_1(U,\bX,\bZ)\}^T]$.
It is easy to show that $E\{ \bolds{\delta}_1(U,\bX,\bZ)\bZ^T \} =0$
and we then have
\begin{eqnarray*}
&&
\bDel_0=E[E(\bZ\bZ^T|U)
\\
&&\hphantom{\bDel_0=E[}{}\times
\{E(\bZ\bZ^T|U)-E(\bZ\bX^T|U)E(\bX\bX^T|U)^{-1}E(\bX\bZ^T|U)\}^{-1}E(\bZ\bZ^T|U)].
\end{eqnarray*}
Therefore,
%e3.5 ###
\begin{equation}
\sqrt{n}(\hat\bbe-\bbe_0)\wcon N(0,R_2(q)\bS^{-1}\bDel_0\bS^{-1}).
\end{equation}

If we replace $R_2(q)$ with $1$ in equations (23) and (24), we end up
with the asymptotic normal distributions of the semi-LS estimators, as
studied in Li and Liang \cite{LiLiang2008}. Thus, $R_2(q)$ determines
the asymptotic relative efficiency (ARE) of semi-CQR relative to
semi-LS. By direct calculations, we see that the ARE for estimating
$\bal(u)$ is $R_2(q)^{-4/5}$ and the ARE for estimating $\bbe$ is
$R_2(q)^{-1}$. It is interesting to see that the same factor, $R_2(q)$,
also appears in the asymptotic efficiency analysis of parametric CQR
\cite{ZouYuan2008} and nonparametric local CQR smoothing
\cite{KaiLiZou2009a}. The basic message is that, with a relatively
large $q$ ($q\ge 9$), $R_2(q)$ is very close to~1 for the normal
errors, but can be much smaller than 1, meaning a huge gain in
efficiency, for the commonly seen non-normal errors. It is also shown
in \cite{KaiLiZou2009a} that $\lim_{q\rightarrow\infty}R_2(q)^{-1}\geq 0.864$ and hence
$\lim_{q\rightarrow\infty}R_2(q)^{-4/5}\geq 0.8896$, which implies that when a large $q$ is used, the ARE is at
least $88.9\%$ for estimating varying-coefficient functions and at
least $86.4\%$ for estimating parametric components.

\begin{remark}
The baseline function estimator $\hat\alp_{0}(u)$ converges to
$\alp_{0}(u)$ plus the average of uniform quantiles of the error
distribution. Therefore, the bias term is zero when the error
distribution is symmetric. Even for nonasymmetric distributions, the
additional bias term converge to the mean of the error, which is zero
for a large value of $q$. Nevertheless, its asymptotic variance differs
from that of the semi-LS estimator by a factor of $R_1(q)$. The study
in Kai, Li and Zou \cite{KaiLiZou2009a} shows that $R_1(q)$ approaches
1 as $q$ becomes large and $R_1(q)$ could be much smaller than 1 with a
smaller $q$ ($q\le 9$) for commonly used non-normal distributions.
\end{remark}

\begin{remark}
The factors $R_1(q)$ and $R_2(q)$ only depend on the error
distribution. We have observed from our simulation study that, as a
function of $q$, the maximum of $R_2(q)$ is often closely approximated
by $R_2$ $(q=9)$. Hence, if we only care about the inference of
$\bal(u)$ and $\bbe$, then $q=9$ seems to be a good default value. On
the other hand, $R_1$ $(q=5)$ is often close to the maximum of $R_1(q)$
based on our numerical study and hence $q=5$ is a good default value
for estimating the baseline function. If prediction accuracy is the
primary interest, then we should use a proper $q$ to maximize the total
contributions from $R_1(q)$ and $R_2(q)$. Practically speaking, one can
choose a $q$ from the interval $[5,9]$ by some popular tuning methods
such as $K$-fold cross-validation. However, we do not expect these CQR
models to have significant differences in terms of model fitting and
prediction because, in many cases, $R_1(q)$ and $R_2(q)$ vary little in
the interval $[5,9]$.
\end{remark}

%s4 ###
\section{Variable selection}\label{s4}

Variable selection is a crucial step in high-dimensional modeling.
Various powerful penalization methods have been developed for variable
selection in parametric models; see Fan and Li \cite{FanLi2006} for a
good review. In the literature, there are only a few papers on variable
selection in semiparametric regression models. Li and Liang
\cite{LiLiang2008} proposed the nonconcave penalized quasi-likelihood
method for variable selection in semiparametric varying-coefficient
models. In this section, we study the penalized semiparametric CQR
estimator.

Let $p_{\lam_n}(\cdot)$ be a pre-specified penalty function with
regularization parameter $\lam_n$. We consider the penalized CQR loss
%e4.1 ###
\begin{equation}\label{ell4cqr}
\sumq\sumn\rho_{\tau_k}\{Y_i-\tilde{\alpha}_{0k}(U_i)-\bX_i^T\tilde{\bal}(U_i)-\bZ_i^T\bbe\}+nq\sum_{j=1}^{d_2}p_{\lam_n}(|\beta_j|).
\end{equation}
By minimizing the above objective function with a proper penalty
parameter $\lam_n$, we can get a sparse estimator of $\bbe$ and hence
conduct variable selection.

Fan and Li \cite{FanLi2001} suggested using a concave penalty function
since it is able to produce an oracular estimator, that is, the
penalized estimator performs as well as if the subset model were known
in advance. However, optimizing (\ref{ell4cqr}) with a concave penalty
function is very challenging because the objective function is
nonconvex and both loss and penalty parts are nondifferentiable.
Various numerical algorithms have been proposed to address the
computational difficulties. Fan and Li \cite{FanLi2001} suggested using
local quadratic approximation (LQA) to substitute for the penalty
function and then optimizing using the Newton--Raphson algorithm.
Hunter and Li \cite{HunterLi2005} further proposed a perturbed version
of LQA to alleviate one drawback of LQA. Recently, Zou and Li
\cite{ZouLi2008} proposed a new unified algorithm based on local linear
approximation (LLA) and further suggested using the one-step LLA
estimator because the one-step LLA automatically adopts a sparse
representation and is as efficient as the fully iterated LLA estimator.
Thus, the one-step LLA estimator is computationally and statistically
efficient.

We proposed to follow the one-step sparse estimate scheme in Zou and Li
\cite{ZouLi2008} to derive a one-step sparse semi-CQR estimator, as
follows. First, we compute the unpenalized semi-CQR estimate
$\hat\bbe{}^{(0)}$, as described in Section \ref{s3}. We then
define
\[%\label{ell5cqr}
G_n(\bbe)=
\sumq\sumn\rho_{\tau_k}\{Y_i-\tilde{\alpha}_{0k}(U_i)-\bX_i^T\tilde{\bal}(U_i)-\bZ_i^T\bbe\}
+nq\sumdd p'_{\lam_n}\bigl(\bigl|\beta_j^{(0)}\bigr|\bigr)|\beta_j|.
\]
We define $\hat{\bbe}{}^{\mathrm{OSE}}=\argmin_{\bbe}G_n(\bbe)$ and call this
the \textit{one-step sparse semi-CQR} estimator. Indeed, this is a
weighted $L_1$ regularization procedure.

We now show that the one-step sparse semi-CQR  estimator $\hat{\bbe}{}^{\mathrm{OSE}}$ enjoys the oracle property. This property holds for a wide
class of concave penalties. To establish the idea, we focus on the SCAD
penalty from Fan and Li \cite{FanLi2001}, which is perhaps the most
popular concave penalty in the literature. Let $\bbe_0=(\bbe_{10}^T,\bbe_{20}^T)^T$ denote the true value of $\bbe$, where
$\bbe_{10}$ is a $s$-vector. Without loss of generality, we assume that
$\bbe_{20}=\bzr$ and that $\bbe_{10}$ contains all nonzero components
of $\bbe_0$. Furthermore, let $\bZ_1$ be the first $s$ elements of
$\bZ$ and define
\[
\bolds{\lambda}(u,\bx,\bz)
=
E[\bZ_1(\bc^T,c\bX^T,\bzr)|U=u]\bD_1^{-1}(u)(I_q,\bone^T\bx,\bone^T\bz)^T.
\]

\begin{thmm}[(Oracle property)]\label{ThmCqrOracle}
Let $p_{\lam}(\cdot)$ be the SCAD penalty. Assume that the
regularity\vspace*{1pt}
conditions \textup{(B1)}--\textup{(B6)} given in the Appendix hold. If $\sqrt{n}\lam_{n}\rar\infty$,
$\lam_{n}\rar 0$,  $nh^4\rightarrow 0$ and
$nh^2/\log(1/h)\rightarrow\infty$ as $n\rightarrow\infty$, then the
one-step semi-CQR estimator $\hat{\bbe}{}^{\mathrm{OSE}}$ must satisfy:
\begin{longlist}
\item[(a)] sparsity, that is, $\hat\bbe{}_2^{\mathrm{OSE}}=0$, with probability tending to one;
\item[(b)] asymptotic normality, that is,
%e4.2 ###
\begin{equation}
\sqrt{n}(\hat\bbe{}_1^{\mathrm{OSE}}-\bbe_{10})
\wcon
N\biggl(0,\frac{1}{c^2}\bS_1^{-1}\bolds{\Lambda}\bS_1^{-1}\biggr),
\end{equation}
where $\bS_1=E(\bZ_1\bZ_1^T)$ and
$\bolds{\Lambda}=\sumq\sumqq\tau_{kk'}E[\{\bZ_1-\bolds{\lambda}_{k}(U,\bX,\bZ)\}\{\bZ_1-\bolds{\lambda}_{k'}(U,\bX,\bZ)\}^T]$ with
$\bolds{\lambda}_{k}(u,\bx,\bz)$ being the $k$th column of the matrix
$\bolds{\lambda}(u,\bx,\bz)$.
\end{longlist}
\end{thmm}

Theorem~\ref{ThmCqrOracle} shows the asymptotic magnitude of the
optimal $\lam_n$. For a given data set with finite sample, it is
practically important to have a data-driven method to select a good
$\lam_n$. Various techniques have been proposed in previous studies,
such as the generalized cross-validation selector \cite{FanLi2001} and
the BIC selector \cite{WangLiTsai2007}. In this work, we use a BIC-like
criterion to select the penalization parameter. The BIC criterion is
defined as
\begin{eqnarray*}
\operatorname{BIC}(\lam)
&=&\log\Biggl(\sumq\sumn\rho_{\tau_k}\{Y_i-\hat\alp_{0k}(U_i)-\bX_i^T\hat\bal(U_i)-\bZ_i^T\hat{\bbe}{}^{\mathrm{OSE}}(\lam)\}\Biggr)
\\
&&{}+
\frac{\log(n)}{n}\mathit{df}_{\lam},
\end{eqnarray*}
where $\mathit{df}_{\lam}$ is the number of nonzero coefficients in the
parametric part of the fitted model. We let $\hat \lam_{\mathrm{BIC}}=\argmin \operatorname{BIC}(\lam)$. The performance of $\hat \lam_{\mathrm{BIC}}$ will be examined in
our simulation studies in the next section.

\begin{remark}
Variable selection in linear quantile regression has been considered in
several papers; see Li and Zhu \cite{Zhu06} and Wu and Liu
\cite{scadqr}. The developed method for sparse semiparametric CQR can
be easily adopted for variable selection in semiparametric quantile
regression. Consider the penalized check loss
%e4.3 ###
\begin{equation}\label{ell4qr}
\sumn\rho_{\tau}\{Y_i-\tilde\alp_{0,\tau}(U_i)-\bX_i^T \tilde\bal_\tau(U_i)-\bZ_i^T\bbe\}+n\sum_{j=1}^{d_2}p_{\lam_n}(|\beta_j|).
\end{equation}
For its one-step version, we use
%e4.4 ###
\begin{equation}\label{ell5qr}
\sumn\rho_{\tau}\{Y_i-\tilde\alp_{0,\tau}(U_i)-\bX_i^T\tilde\bal_\tau(U_i)-\bZ_i^T\bbe\}+n\sumdd p'_{\lam_n}\bigl(\bigl|\beta_j^{(0)}\bigr|\bigr)|\beta_j|,
\end{equation}
where $\hat\bbe{}^{(0)}$ denotes the unpenalized semiparametric quantile
regression estimator defined in Section \ref{s2}. We can also prove the
oracle property of the one-step sparse semiparametric quantile
regression estimator by following the lines of proof for
Theorem~\ref{ThmCqrOracle}. For reasons of brevity, we omit the details
here.
\end{remark}

%s5 ###
\section{Numerical studies}\label{s5}

In this section, we conduct simulation studies to assess the
finite-sample performance of the proposed procedures and illustrate the
proposed methodology on a real-world data set in a health study. In all
examples, we fix the kernel function to be the Epanechnikov kernel,
that is, $K(u)=\frac{3}{4}(1-u^2)_{+}$, and we use the SCAD penalty
function for variable selection. Note that all proposed estimators,
including semi-QR, semi-CQR and one-step sparse semi-CQR, can be
formulated as linear programming (LP) problems. In our study, we solved
these estimators by using LP tools.

\begin{example}\label{ex41}
In this example, we generate 400 random samples, each consisting of
$n=200$ observations, from the  model
\[
Y=\alpha_1(U)X_1+\alpha_2(U)X_2+\beta_1 Z_1+\beta_2 Z_2+\beta_3 Z_3+\eps,
\]
where $\alpha_1(U)=\sin(6\pi U)$, $\alpha_2(U)=\sin(2\pi U)$,
$\beta_1=2$, $\beta_2=1$ and $\beta_3=0.5$. The covariate $U$ is from
the uniform distribution on $[0,1]$. The covariates $X_1,X_2,Z_1,Z_2$
are jointly normally distributed with mean 0, variance 1 and
correlation $2/3$. The covariate $Z_3$ is Bernoulli with
$\Pr(Z_3=1)=0.4$. Furthermore, $U$ and $(X_1,X_2,Z_1,Z_2,Z_3)$ are
independent. In our simulation, we considered the following error
distributions: $N(0,1)$, logistic, standard Cauchy, $t$-distribution
with 3 degrees of freedom, mixture of normals $0.9N(0,1)+0.1 N(0,10^2)$ and log-normal distribution. Because the error is
independent of the covariates, the least-squares (LS), quantile
regression (QR) and composite quantile regression (CQR) procedures
provide estimates for the same quantity and hence are directly
comparable.

\textit{Performance of $\hat{\bbe}_{\tau}$ and $\hat{\bbe}$.}
%We first investigate the effect of bandwidth choice.
To examine the performance of the proposed procedures with a wide range
of bandwidths, three bandwidths for LS were considered, $h=0.085$,
$0.128$, $0.192$, which correspond to the undersmoothing, appropriate
smoothing and oversmoothing, respectively. By straightforward
calculation, as in Kai, Li and Zou \cite{KaiLiZou2009a}, we can produce
two simple formulas for the asymptotic optimal bandwidths for QR and
CQR: $h_{\mathrm{CQR}}=h_{\LS}\cdot\Rb^{1/5}$ and
$h_{\mathrm{QR},\tau}=h_{\LS}\cdot\{\tau(1-\tau)/f[F^{-1}(\tau)]\}^{1/5}$, where $h_{\LS}$ is the
asymptotic optimal bandwidth for LS. We considered only  the case of
normal error. The bias and standard deviation based on 400 simulations
are reported in Table~\ref{tabMeanSD}. First, we see that the
estimators are not very sensitive to the choice of bandwidth. As for
the estimation accuracy, all three estimators have comparable bias and
the differences are shown in standard deviation. The LS estimates have
the smallest standard deviation, as expected. The CQR estimates are
slightly worse than the LS estimates.

In the second study, we fixed $h=0.128$ and compared the efficiency of
QR and CQR relative to LS. Reported in Table~\ref{tabRMSE} are RMSEs,
the ratios of the MSEs of the QR and CQR estimators to the LS estimator
for different error distributions. Several observations can be made
from Table~\ref{tabRMSE}. When the error follows the normal
distribution, the RMSEs of CQR are slightly less than 1. For all other
non-normal distributions in the table, the RMSE can be much greater
than 1, indicating a huge gain in efficiency. These findings agree with
the asymptotic theory. For QR estimators, their performance varies and
depends heavily on the level of quantile and the error distribution.
Overall, CQR outperforms both QR and LS.

\textit{Performance of $\hat{\bal}_{\tau}$ and $\hat{\bal}$}.
We now compare the LS, QR and CQR estimates for $\bal$ by using the
ratio of average squared errors (RASE). We first compute
\[
\operatorname{ASE}=\Biggl\{\frac 1{\ngrid}\sum_{m=1}^{d_1}\sum_{k=1}^{\ngrid}\{\hat{a}_m(u_k)-a_m(u_k)\}^2\Biggr\},
\]
where $\{u_k\dvtx k=1,\ldots,\ngrid\}$ is a set of grid points uniformly
placed on $[0,1]$ with $\ngrid=200$. RASE is then defined to be
%e5.1 ###
\begin{equation}
\operatorname{RASE}(\hat{g})
=
\frac{\operatorname{ASE}(\hat{g}_{\LS})}{\operatorname{ASE}(\hat{g})}
\end{equation}
for an estimator $\hat{g}$, where $\hat{g}_{\LS}$ is the
least-squares-based estimator.

The sample mean and standard deviation of the RASEs over 400
simulations are presented in Table~\ref{tabRASE}, where the values in
the parentheses are the standard deviations. The findings are quite
similar to those in Table~\ref{tabRMSE}.%t1 ###
\begin{table}
\caption{Summary of the bias and standard deviation over 400
simulations}\label{tabMeanSD}
\begin{tabular}{@{}lcccc@{}}
\hline
&&\multicolumn{3}{c@{}}{$\bolds{\operatorname{Bias}}\bolds{(}\mathbf{SD}\bolds{)}$}\\[-5pt]
&&\multicolumn{3}{c@{}}{\hrulefill}\\
$\bolds{h}$& \multicolumn{1}{c}{\textbf{Method}}&\multicolumn{1}{c}{$\bolds{\hat\beta}_\mathbf{1}$}&\multicolumn{1}{c}{$\bolds{\hat\beta}_\mathbf{2}$}&\multicolumn{1}{c@{}}{$\bolds{\hat\beta}_\mathbf{3}$}\\
\hline
0.085&LSE&$-$0.012 (0.128)&\phantom{$-$}0.008 (0.121)&$-$0.009 (0.171)\\
&CQR$_9$&$-$0.009 (0.131)&\phantom{$-$}0.009 (0.125)&$-$0.007 (0.172)\\
&QR$_{0.25}$&$-$0.017 (0.163)&\phantom{$-$}0.009 (0.161)&$-$0.151 (0.237)\\
&QR$_{0.50}$&$-$0.012 (0.155)&\phantom{$-$}0.011 (0.151)&$-$0.007 (0.198)\\
&QR$_{0.75}$&$-$0.007 (0.165)&\phantom{$-$}0.005 (0.158)&\phantom{$-$}0.122 (0.216)\\[3pt]
0.128&LSE&$-$0.009 (0.121)&\phantom{$-$}0.005 (0.117)&$-$0.008 (0.164)\\
&CQR$_9$&$-$0.010 (0.127)&\phantom{$-$}0.008 (0.121)&$-$0.005 (0.163)\\
&QR$_{0.25}$&$-$0.010 (0.159)&\phantom{$-$}0.003 (0.152)&$-$0.082 (0.227)\\
&QR$_{0.50}$&$-$0.008 (0.154)&\phantom{$-$}0.011 (0.147)&$-$0.004 (0.193)\\
&QR$_{0.75}$&$-$0.012 (0.163)&\phantom{$-$}0.003 (0.161)&\phantom{$-$}0.071 (0.207)\\[3pt]
0.192&LSE&$-$0.007 (0.128)&\phantom{$-$}0.001 (0.123)&$-$0.008 (0.169)\\
&CQR$_9$&$-$0.009 (0.131)&\phantom{$-$}0.005 (0.127)&$-$0.005 (0.169)\\
&QR$_{0.25}$&$-$0.006 (0.169)&$-$0.004 (0.169)&$-$0.061 (0.230)\\
&QR$_{0.50}$&$-$0.005 (0.153)&\phantom{$-$}0.006 (0.152)&$-$0.007 (0.191)\\
&QR$_{0.75}$&$-$0.012 (0.170)&\phantom{$-$}0.007 (0.171)&\phantom{$-$}0.049 (0.225)\\
\hline
\end{tabular}
\end{table}
We see that CQR performs
almost as well as LS when the error is normally distributed. Also, its
RASEs are much larger than 1 for other non-normal error distributions.
%t2 ###
\begin{table}
\caption{Summary of the ratio of MSE over 400 simulations}\label{tabRMSE}
\begin{tabular*}{190pt}{@{\extracolsep{\fill}}lccc@{}}
\hline
&\multicolumn{3}{c@{}}{$\mathbf{RMSE}$}\\[-5pt]
&\multicolumn{3}{c@{}}{\hrulefill}\\
\multicolumn{1}{@{}l}{$\mathbf{Method}$}&\multicolumn{1}{c}{$\bolds{\hat\beta}_\mathbf{1}$}&\multicolumn{1}{c}{$\bolds{\hat\beta}_\mathbf{2}$}&\multicolumn{1}{c@{}}{$\bolds{\hat\beta}_\mathbf{3}$}\\
\hline
\multicolumn{4}{@{}l@{}}{Standard normal}\\
\quad CQR$_9$&0.920&0.932&1.011\\
\quad QR$_{0.25}$&0.585&0.594&0.460\\
\quad QR$_{0.50}$&0.621&0.631&0.724\\
\quad QR$_{0.75}$&0.554&0.528&0.561\\[3pt]
\multicolumn{4}{@{}l@{}}{Logistic}\\
\quad CQR$_9$&1.044&1.083&1.016\\
\quad QR$_{0.25}$&0.651&0.664&0.502\\
\quad QR$_{0.50}$&0.826&0.871&0.799\\
\quad QR$_{0.75}$&0.661&0.732&0.527\\[3pt]
\multicolumn{4}{@{}l@{}}{Standard Cauchy}\\
\quad CQR$_9$&15,246&106,710&52,544\\
\quad QR$_{0.25}$&\phantom{0,}8894&\phantom{0}56,704&24,359\\
\quad QR$_{0.50}$&19,556&137,109&66,560\\
\quad QR$_{0.75}$&\phantom{0,}8223&\phantom{0}62,282&26,210\\[3pt]
\multicolumn{4}{@{}l@{}}{$t$-distribution with $\mathit{df}=3$}\\
\quad CQR$_9$&1.554&1.546&1.683\\
\quad QR$_{0.25}$&1.000&0.948&0.819\\
\quad QR$_{0.50}$&1.354&1.333&1.451\\
\quad QR$_{0.75}$&0.935&1.059&0.859\\[3pt]
\multicolumn{4}{@{}l@{}}{$0.9 N(0,1)+0.1 N(0,10^2)$}\\
\quad CQR$_9$&5.752&4.860&5.152\\
\quad QR$_{0.25}$&3.239&3.096&2.300\\
\quad QR$_{0.50}$&5.430&4.730&4.994\\
\quad QR$_{0.75}$&3.790&2.952&2.515\\[3pt]
\multicolumn{4}{@{}l@{}}{Log-normal}\\
\quad CQR$_9$&3.079&3.369&3.732\\
\quad QR$_{0.25}$&5.198&5.361&3.006\\
\quad QR$_{0.50}$&2.787&2.829&3.139\\
\quad QR$_{0.75}$&0.819&0.868&0.823\\
\hline
\end{tabular*}
\end{table}
The efficiency gain can be substantial. Note that for the Cauchy
distribution, RASEs of QR and CQR are huge---this  is because LS fails
when the error variance is infinite.
\end{example}

\begin{example}\label{42}
The goal is to compare the proposed one-step sparse semi-CQR estimator
with the one-step sparse semi-LS estimator.
%t3 ###
\begin{table}
\tabcolsep=0pt
\caption{Summary of the RASE over 400 simulations}\label{tabRASE}
\begin{tabular*}{\textwidth}{@{\extracolsep{\fill}}lcccccc@{}}
\hline
&\textbf{Normal}&\textbf{Logistic}&\textbf{Cauchy}&{$\bolds{t}_{\mathbf{3}}$}&\textbf{Mixture}&\textbf{Log-normal}\\
\hline
CQR$_9$     &0.968 (0.104)&1.040 (0.134)&$12{,}872$ (176719)&1.428 (1.299)&3.292 (1.405)&2.455 (1.498)\\
QR$_{0.25}$ &0.666 (0.160)&0.720 (0.203)&\phantom{0,}7621 (110692)&0.958 (0.647)&2.029 (1.003)&3.490 (3.224)\\
QR$_{0.50}$ &0.771 (0.184)&0.881 (0.206)&$13{,}720$ (187298)&1.274 (1.166)&3.155 (1.323)&2.155 (1.674)\\
QR$_{0.75}$ &0.681 (0.191)&0.713 (0.201)&\phantom{0}5781 (87909)&0.896 (0.325)&1.953 (0.905)&0.824 (0.679)\\
\hline
\end{tabular*}
\end{table}
In this example, 400 random
samples, each consisting of $n=200$ observations, were generated from
the varying-coefficient partially linear model
\[
Y=\alpha_1(U)X_1+\alpha_2(U)X_2+\bbe^T\bZ+\eps,
\]
where $\bbe=(3,1.5,0,0,2,0,0,0)^T$ and the covariate vector $(X_1, X_2,\bZ^T)^T$ is normally distributed with mean 0, variance 1 and
correlation $0.5^{|i-j|}$ $(i,j=1,\ldots,10)$. Other model settings are
exactly the same as those in Example~\ref{ex41}. We use the generalized
mean square error (GMSE), as defined in
\cite{LiLiang2008},
%e5.2 ###
\begin{equation}
\operatorname{GMSE}(\hat\bbe)=(\hat\bbe-\bbe)^T E(\bZ\bZ^T)(\hat\bbe-\bbe),
\end{equation}
to assess the performance of variable selection procedures for the
parametric component. For each procedure, we calculate the relative
GMSE (RGMSE), which is defined to be the ratio of the GMSE of a
selected final model to that of the unpenalized least-squares estimate
under the full model.

%t4 ###
\begin{table}[b]
\caption{One-step estimates for variable selection in semiparametric
models}\label{Variableselection}
\begin{tabular*}{\textwidth}{@{\extracolsep{\fill}}lcccccc@{}}
\hline
&\multicolumn{1}{c}{\multirow{2}{33pt}{\textbf{RGMSE}}}&\multicolumn{2}{c}{\textbf{No.~of zeros}}&\multicolumn{3}{c@{}}{\textbf{Proportion of fits}}\\[-5pt]
&&\multicolumn{2}{c}{\hrulefill}&\multicolumn{3}{c@{}}{\hrulefill}\\
\textbf{Method}&\textbf{Median (MAD)}&\textbf{C}&\textbf{IC}&\textbf{U-fit}&\textbf{C-fit}&\textbf{O-fit}\\
\hline
\multicolumn{4}{@{}l@{}}{Standard normal}\\
\quad One-step LS &0.335 (0.194)&4.825&0.000&0.000&0.867&0.133\\
\quad One-step CQR &0.288 (0.213)&4.990&0.000&0.000&0.990&0.010\\[3pt]
\multicolumn{4}{@{}l@{}}{Logistic}\\
\quad One-step LS&0.352 (0.197)&4.805&0.000&0.000&0.870&0.130\\
\quad One-step CQR&0.289 (0.206)&4.975&0.000&0.000&0.975&0.025\\[3pt]
\multicolumn{4}{@{}l@{}}{Standard Cauchy}\\
\quad One-step LS&0.956 (0.249)&2.920&0.795&0.595&0.108&0.297\\
\quad One-step CQR&0.005 (0.021)&5.000&0.295&0.210&0.790&0.000\\[3pt]
\multicolumn{4}{@{}l@{}}{$t$-distribution with $\mathit{df}=3$} \\
\quad One-step LS&0.346 (0.179)&4.803&0.000&0.000&0.860&0.140\\
\quad One-step CQR&0.183 (0.177)&4.987&0.000&0.000&0.988&0.013\\[3pt]
\multicolumn{4}{@{}l@{}}{$0.9 N(0,1) + 0.1 N(0,10^2)$}\\
\quad One-step LS&0.331 (0.190)&4.848&0.000&0.000&0.883&0.117\\
\quad One-step CQR&0.060 (0.083)&4.997&0.000&0.000&0.998&0.003\\[3pt]
\multicolumn{4}{@{}l@{}}{Log-normal}\\
\quad One-step LS&0.303 (0.182)&4.845&0.000&0.000&0.887&0.113\\
\quad One-step CQR&0.111 (0.118)&4.990&0.000&0.000&0.990&0.010\\
\hline
\end{tabular*}
\end{table}

The results over 400 simulations are summarized in
Table~\ref{Variableselection}, where the column ``RGMSE'' reports both
the median and MAD of 400 RGMSEs. Both columns ``C'' and ``IC'' are
measures of model complexity. Column ``C'' shows the average number of
zero coefficients correctly estimated to be zero and column ``IC''
presents the average number of nonzero coefficients incorrectly
estimated to be zero. In the column labeled ``U-fit'' (short for
``under-fit''), we present the proportion of trials excluding any
nonzero coefficients in 400 replications. Likewise, we report the
probability of trials selecting the exact subset model and the
probability of trials including all three significant variables and
some noise variables in the columns ``C-fit'' (``correct-fit'') and
``O-fit'' (``over-fit''), respectively. From
Table~\ref{Variableselection}, we see that both variable selection
procedures dramatically reduce model errors, which clearly show the
virtue of variable selection. Second, the one-step CQR performs better
than the one-step LS in terms of all of the criteria: RGMSE, number of
zeros and proportion of fits, and for all of the error distributions in
Table~\ref{Variableselection}. It is also interesting to see that in
the normal error case, the one-step CQR seems to perform no worse than
the one-step LS (or even slightly better). We performed the
Mann--Whitney test to compare their RGMSEs and the corresponding
$p$-value is $0.0495$. This observation appears to be contradictory to
the asymptotic theory. However, this ``contradiction'' can be explained
by observing that the one-step CQR has better variable selection
performance. Note that the one-step CQR has significantly higher
probability of correct-select than the one-step LS, which also tends to
overselect. Thus, the one-step LS needs to estimate a larger model than
the truth, compared to the one-step CQR.
\end{example}

\begin{example}\label{43}
As an illustration, we apply the proposed procedures to
analyze the plasma
beta-carotene level data set %studied in Wang and Li \citet{WangLi2009}.
collected by a cross-sectional study \cite{NierenbergEtal1989}. This
data set  consists of 273 samples. Of interest are  the relationships
between the plasma beta-carotene level and the following  covariates:
age, smoking status, quetelet index (BMI), vitamin use, number of calories, grams
of fat, grams of fiber, number of alcoholic drinks, cholesterol and
dietary beta-carotene. The complete description of the data can be
found in the StatLib database via the link \href{http://lib.stat.cmu.edu/datasets/Plasma\_Retinol}{lib.stat.cmu.edu/datasets/Plasma\_Retinol}.

We fit the data by using a partially linear model with $U$ being
``\textit{dietary beta-carotene}.'' The covariates ``\textit{smoking
status}'' and ``\textit{vitamin use}'' are categorical and are thus
replaced with dummy variables. All of the other covariates are
standardized. We applied the one-step sparse CQR and LS estimators to
fit the partially linear regression model. Five-fold cross-validation
was used to select the bandwidths for LS and CQR. We used the first 200
observations as a training data set to fit the model and to select
significant variables, then used the remaining 73 observations to
evaluate the predictive ability of the selected model.

%t5 ###
\begin{table}
\caption{Selected parametric components for plasma beta-carotene level
data}\label{ex3tb1}
\begin{tabular*}{250pt}{@{\extracolsep{\fill}}ld{4.2}d{2.2}@{}}
\hline
&\multicolumn{1}{c}{${\bolds{\hat\beta}}{}^{\mathbf{OSE}}_{\mathbf{LS}}$}&\multicolumn{1}{c@{}}{${\bolds{\hat\beta}}{}^{\mathbf{OSE}}_{\mathbf{CQR}}$}\\
\hline
Age&0&0\\
Quetelet index&0&0\\
Calories&-100.47&0\\
Fat&52.60&0\\
Fiber&87.51&29.89\\
Alcohol&44.61&0\\
Cholesterol&0&0\\
Smoking status (never)&51.71&0\\
Smoking status (former)&72.48&0\\
Vitamin use (yes, fairly often)&130.39&30.21\\
Vitamin use (yes, not often)&0&0\\[3pt]
MAPE&111.28&58.11\\
\hline
\end{tabular*}
\end{table}

The prediction performance is measured by the median absolute
prediction error (MAPE), which is the median of $\{|y_i - \hat y_i|,
i=1,2,\ldots,73\}$. To see the effect of $q$ on the CQR estimate, we tried
$q=5,7,9$. We found that the selected $Z$-variables are the same for
these three values of $q$ and their MAPEs are $58.52$, $58.11$,
$62.43$, respectively. Thus, the effect of $q$ is minor. The resulting
model with $q=7$ is given in Table~\ref{ex3tb1} and the estimated
intercept function is depicted in Figure~\ref{ex3fg1}.
%f1 ###
\begin{figure}[b]

\includegraphics{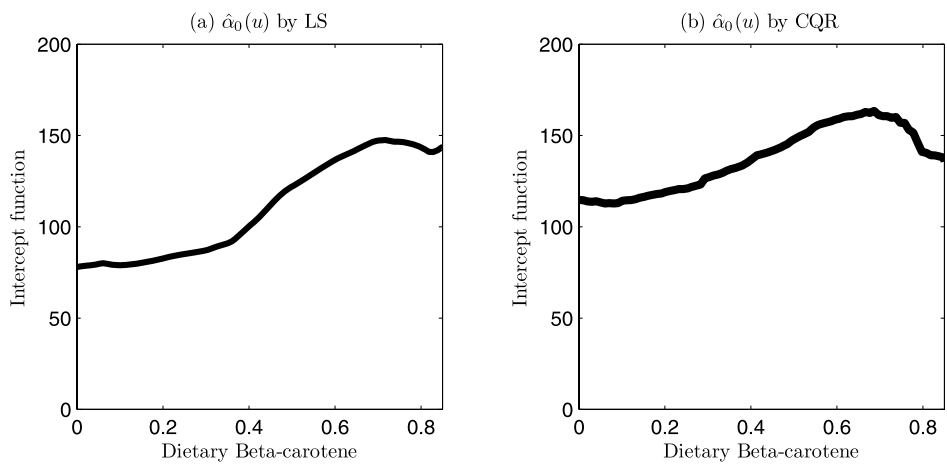}
 \caption{Plot of estimated intercept function of
\textit{dietary beta-carotene}: \textup{(a)} the estimated intercept
function by LS method; \textup{(b)}  the estimated intercept function
by CQR method with $q=7$.} \label{ex3fg1}
\end{figure}
From
Table~\ref{ex3tb1}, it can be seen that the CQR model is much sparser
than the LS model. Only two covariates, ``\textit{fiber consumption per
day}'' and ``\textit{fairly often use of vitamin}'' are included in the
parametric part of the CQR model. Meanwhile, the CQR model has much
better prediction performance than the LS model, whose MAPE is
$111.28$.
\end{example}

%s6 ###
\section{Discussion}\label{s6}

We discuss some directions in which this work could be further
extended. We have focused on using uniform weights in composite
quantile regression. In theory, we can use nonuniform weights, which
may provide an even more efficient estimator when a reliable estimate
of the error distribution is available. Koenker \cite{Koenker1984}
discussed the theoretically optimal weights. Bradic, Fan and Wang
\cite{Fan2010} suggested a data-driven weighted CQR for parametric
linear regression, in which the weights mimic the optimal weights. The
idea in Bradic, Fan and Wang~\cite{Fan2010} can be easily extended to
the semi-CQR estimator, which will be investigated in detail in a
future paper.

Penalized Wilcoxon rank regression has been considered independently in Leng \cite{Leng10} and
Wang and Li \cite{WangLi2009} and found to
achieve a similar efficiency property of CQR for variable selection in
parametric linear regression. We could also generalize rank regression
to handle semiparametric varying-coefficient partially linear models.
In a working paper, we show that rank regression is exactly equivalent
to CQR using $q=n-1$ quantiles with uniform weights. This result
indicates that CQR is more flexible than rank regression because we can
easily use flexible nonuniform weights in CQR to further improve
efficiency, as in Bradic, Fan and Wang \cite{Fan2010}. Obviously, CQR
is also computationally more efficient than rank regression. We note
that in parametric linear regression models, rank regression has no
efficiency gain over least-squares for estimating the intercept. This
result is expected to hold for estimating the baseline function in the
semiparametric varying-coefficient partially linear model.

When the number of varying coefficient components is large, it is also
desirable to consider selecting a few important components. This
problem was studied in Wang and Xia \cite{WangXia09}, where a
LASSO-type penalized local least-squares estimator was proposed. It
would be interesting to apply CQR to their method to further improve
the estimation efficiency.

%s7 ###
\section{Proofs}\label{s7}

To establish the asymptotic properties of the proposed estimators, the
following regularity conditions are imposed:
\begin{longlist}
\item[(C1)] the random variable $U$ has  bounded support $\Omega$ and its density function $f_U(\cdot)$ is positive and has a continuous second derivative;
\item[(C2)] the varying coefficients $\alp_0(\cdot)$ and $\bal(\cdot)$ have continuous second derivatives in $u \in \Omega$;
\item[(C3)] $K(\cdot)$ is a symmetric density function with bounded support and satisfies a Lipschitz condition;
\item[(C4)] the random vector $\bZ$ has bounded support;  \item[(C5)] for the semi-QR procedure,
\begin{longlist}[(C4)]
\item[(i)] $F_{\tau}(0|u,\bx,\bz)=\tau$ for all $(u,\bx,\bz),$ and $f_{\tau}(\cdot|u,\bx,\bz)$ is bounded away from zero and has a continuous and uniformly bounded derivative;
\item[(ii)] $\bA_1(u)$ defined in Theorem~\ref{ThmQrAlpha1} and $\bA_2(u)$ defined in Theorem~\ref{ThmQrAlpha2} are nonsingular for all $u\in\Omega$;
\end{longlist}
\item[(C6)] for the semi-CQR procedure,
\begin{longlist}[(C4)]
\item[(i)] $f(\cdot)$ is bounded away from zero and has a continuous and uniformly bounded derivative;
\item[(ii)] $\bD_1(u)$ defined in Theorem~\ref{ThmCqrAlpha1} and $\bD_2(u)$ defined in Theorem~\ref{ThmCqrAlpha}  are nonsingular for all $u\in\Omega$.
\end{longlist}
\end{longlist}

Although the proposed semi-QR and semi-CQR procedures require
different regularity conditions, the proofs follow similar strategies.
For brevity, we only present the detailed proofs for the semi-CQR
procedure. The detailed proofs for the semi-QR procedure was given in
the earlier version of this paper. Lemma~\ref{lemma1} below, which is a
direct result of Mack and Silverman \cite{MarkSilverman1982}, will be
used repeatedly  in our proofs. Throughout the proofs, identities of
the form $G(u)=O_p(a_n)$ always stand for $\sup_{u\in\Omega}
|G(u)|=O_p(a_n)$.

\begin{lemma}\label{lemma1}
Let $(\bX_1,Y_1),\ldots,(\bX_n,Y_n)$ be i.i.d.~random vectors, where the $Y_i$'s
are scalar random variables. Assume, further, that $E|Y|^r<\infty$ and
that $\sup_{\bx}\int|y|^r f(\bx,y)\,dy<\infty$, where $f$ denotes the
joint density of $(\bX,Y)$. Let $K$ be a bounded positive function with
bounded support, satisfying a Lipschitz condition. Then,
\[
\sup_{\bx\in D}\Biggl|n^{-1}\sumn\{K_h(\bX_i-\bx)Y_i-E[K_h(\bX_i-\bx)Y_i]\}\Biggr|=O_p\biggl(\frac{\log^{1/2}(1/h)}{\sqrt{nh}}\biggr),
\]
provided that $n^{2\eps-1}h\rightarrow\infty$ for some $\eps<1-r^{-1}$.
\end{lemma}

Let $\eta_{i,k}=I(\eps_i\leq c_k)-\tau_k$ and $\eta^*_{i,k}(u)=I\{\eps_i\leq c_k-r_i(u)\}-\tau_k$,
where $r_{i}(u)=\alp_0(U_i)-\alp_0(u)-\alp'_0(u)(U_i-u)+\bX_i^T\{\bal(U_i)-\bal(u)-\bal'(u)(U_i-u)\}$.
Furthermore, let
$\tilde\bth{}^*=\sqrt{nh}\{\tilde a_{01}-\alp_0(u)-c_1,\ldots,
\tilde a_{0q}-\alp_0(u)-c_q,\{\tilde\ba-\bal(u)\}^T,\{\tilde \bbe-\bbe_0\}^T,h\{\tilde b_0-\alp'_0(u)\},h\{\tilde\bb-\bal'(u)\}^T\}^T$
and $\bX^*_{i,k}(u)=\{\be_k^T,\bX_i^T,\break \bZ_i^T,(U_{i}-u)/h,\bX_i^T(U_{i}-u)/h \}^T$,
where $\be_{k}$ is a $q$-vector with 1 at the $k$th position and 0 elsewhere.

In the proof of Theorem \ref{ThmCqrAlpha1}, we will first show the
following asymptotic representation of $\{\tilde\ba_0,\tilde b_0,\tilde\ba,\tilde\bb,\tilde\bbe\}$:
%e7.1 ###
\begin{equation}
\tilde\bth{}^*=-f^{-1}_U(u)\{\bS^*(u)\}^{-1}\bW^*_n(u)+O_p\bigl(h^2+\log^{1/2}(1/h)/\sqrt{nh}\bigr),
\end{equation}
where $\bS^*(u)=\diag\{\bD_1(u),c\mu_2\bB_2(u)\}$ and
\[
\bW^*_n(u)=\frac{1}{\sqrt{nh}}\sumq\sumn K\{(U_i-u)/h\}\eta^*_{i,k}(u)\bX_{i,k}^*(u).
\]
The asymptotic normality of $\{\tilde\ba_0,\tilde b_0,\tilde\ba,\tilde\bb,\tilde\bbe\}$ then follows by demonstrating the asymptotic
normality of $\bW^*_n(u)$.

\begin{pf*}{Proof of Theorem \ref{ThmCqrAlpha1}}
Recall that $\{\tilde\ba_0,\tilde\ba,\tilde\bbe,\tilde b_0,\tilde\bb\}$ minimizes
\[
\sumq\sumn\rho_{\tau_k}[Y_i-a_{0k}-b_0(U_i-u)-\bX_i^T\{\ba+\bb(U_i-u)\}-\bZ_i^T\bbe]K_h(U_i-u).
\]
We write $Y_i-a_{0k}-b_0(U_i-u)-\bX_i^T\{\ba+\bb(U_i-u)\}-\bZ_i^T\bbe=(\eps_i-c_k)+r_{i}(u)-\Delta_{i,k}$, where
$\Delta_{i,k}=\{\bX_{i,k}^*(u)\}^T\bth^*/\sqrt{nh}$. Then,
$\tilde\bth{}^*$ is also the minimizer of
\[
L^*_n(\bth^*)=\sumq\sumn K_i(u)[\rho_{\tau_{k}}\{(\eps_i-c_k)+r_i(u)-\Delta_{i,k}\}-\rho_{\tau_{k}}\{(\eps_i-c_k)+r_i(u)\}],
\]
where $K_i(u)=K\{(U_i-u)/h\}$. By applying the\vspace*{-1pt} identity \cite{Knight1998}
%e7.2 ###
\begin{equation}\label{Knight}
\qquad\rho_{\tau}(x-y)-\rho_{\tau}(x)=y\{I(x\leq 0)-\tau\}+\int_0^{y}\{I(x\leq z)-I(x\leq 0)\}\,dz,\vspace*{-1pt}
\end{equation}
we\vspace*{-1pt} have
\begin{eqnarray*}
L^*_n(\bth^*)
&=&
\sumq\sumn K_i(u)\biggl\{\Delta_{i,k}[I\{\eps_i\leq c_k -r_i(u)\}-\tau_k]\vspace*{-1pt}
\\
&&\hphantom{\sumq\sumn K_i(u)\ }{}+
\int_{0}^{\Delta_{i,k}}[I\{\eps_i\leq c_k-r_i(u)+z\}-I\{\eps_i\leq c_k-r_i(u)\}]\,dz\biggr\}
\\
% &=& ( \frac{1}{\sqrt{nh}} \sumq  \sumn K_i(u) \eta^*_{i,k}(u) \bX_{i,k}^*(u)  )^T \bth^* +  \sumq B^*_{n,k}(\bth^*)  \nn\\
&=&\{\bW^*_n(u)\}^T\bth^*+\sumq B^*_{n,k}(\bth^*),\vspace*{-1pt}
\end{eqnarray*}
where
\[
B^*_{n,k}(\bth^*)=\sumn K_i(u)\int_{0}^{\Delta_{i,k}}[I\{\eps_i\leq c_k-r_i(u)+z\}-I\{\eps_i\leq c_k-r_i(u)\}]\,dz.\vspace*{-1pt}
\]
Since $B^*_{n,k}(\bth^*)$ is a\vspace*{1.5pt} summation of i.i.d.~random variables of
the kernel form, it follows, by\vspace*{-1pt} Lemma~\ref{lemma1}, that
\[
B^*_{n,k}(\bth^*)=E[B^*_{n,k}(\bth^*)]+O_p\bigl(\log^{1/2}(1/h)/\sqrt{nh}\bigr).\vspace*{-1pt}
\]
%So
The conditional expectation of $\sumq B^*_{n,k}(\bth^*)$ can be
calculated\vspace*{-1pt} as
\begin{eqnarray*}
%  &&\sumq E[ B^*_{n,k}(\bth^*)]\nn\\
&&
\sumq E[B^*_{n,k}(\bth^*)|U,\bX,\bZ]\vspace*{-1pt}
\\
&&\qquad=
\sumq\sum_{i=1}^{n}K_i(u)\int_{0}^{\Delta_{i,k}}\bigl[F\bigl(c_k-r_i(u)+z\bigr)-F\bigl(c_k-r_i(u)\bigr)\bigr]\,dz
\\
%   &=&\sumq   \sum_{i=1}^{n} K_i(u)
%         \int_{0}^{\Delta_{i,k}} [  z f(c_k - r_i(u)) + O(z^2)
%        ]   dz   \nonumber \\
%   &=&\sumq   \sum_{i=1}^{n} K_i(u)
%         [\Delta_{i,k}^2 f(c_k - r_i(u))/2  + O(\Delta_{i,k}^3)] \nonumber \\
&&\qquad=
\frac{1}{2}(\bth^*)^T\Biggl(\frac{1}{nh}\sumq\sumn K_i(u)f\bigl(c_k-r_i(u)\bigr)\{\bX_{i,k}^*(u)\}\{\bX_{i,k}^*(u)\}^T\Biggr)\bth^*
\\
&&\qquad\quad{}+
O_p\bigl(\log^{1/2}(1/h)/\sqrt{nh}\bigr)
\\
&&\qquad\triangleq
\frac{1}{2}(\bth^*)^T\bS^*_{n}(u)\bth^*+ O_p\bigl(\log^{1/2}(1/h)/\sqrt{nh}\bigr).\vspace*{-1pt}
\end{eqnarray*}
%Define $R^*_{n,k}(\bth^*)= B^*_{n,k}(\bth^*) - E[ B^*_{n,k}(\bth^*)|U,\bX,\bZ] $. It can be shown in Lemma 3 that
%    R^*_{n,k}(\bth^*)=O_p(1/\sqrt[4]{nh}).
Then,\vspace*{-1pt}
\begin{eqnarray*}%\label{Ln4}
L^*_n(\bth^*)
&=&
\{\bW^*_n(u)\}^T\bth^*+\sumq E[B^*_{n,k}(\bth^*)]+O_p\bigl(\log^{1/2}(1/h)/\sqrt{nh}\bigr)
\\
&=&
\{\bW^*_n(u)\}^T\bth^*+\sumq E\{E[B^*_{n,k}(\bth^*)|U,\bX,\bZ]\}+O_p\bigl(\log^{1/2}(1/h)/\sqrt{nh}\bigr)
\\
&=&
\{\bW^*_n(u)\}^T\bth^*+\frac{1}{2}(\bth^*)^T E[\bS^*_{n}(u)]\bth^*+O_p\bigl(\log^{1/2}(1/h)/\sqrt{nh}\bigr).
\end{eqnarray*}
It can be shown that $ E [\bS^*_n(u)] = f_U(u) \bS^*(u) + O(h^2)$.
Therefore, we can write $L_n(\bth^*)$ as
\[
L^*_n(\bth^*)=\{\bW^*_n(u)\}^T\bth^*+\frac{f_U(u)}{2}(\bth^*)^T\bS^*(u)\bth^*+O_p\bigl(h^2+\log^{1/2}(1/h)/\sqrt{nh}\bigr).
\]
By applying the convexity lemma \cite{Pollard1991} and the quadratic
approximation lemma~\cite{FanGijbels1996}, the minimizer of
$L^*_n(\bth^*)$ can be expressed as
%e7.3 ###
\begin{equation}\label{CQRalphaR}
    \tilde \bth{}^* = -f^{-1}_U(u) \{\bS^*(u)\}^{-1} \bW^*_n(u) + O_p\bigl(h^2+\log^{1/2}(1/h)/\sqrt{nh}\bigr),
\end{equation}
which holds uniformly for $u \in \Omega$. Meanwhile, for any point $u
\in \Omega$, we have
%e7.4 ###
\begin{equation}
\tilde\bth{}^*=-f^{-1}_U(u)\{\bS^*(u)\}^{-1}\bW^*_n(u)+o_p(1).
\end{equation}
Note that $\bS^*(u)=\diag \{ \bD_1(u),  c \mu_2 \bB_2(u) \}$ is a
quasi-diagonal matrix. So,
%e7.5 ###
\begin{equation}\label{CQRalphaR1}
\sqrt{nh}
\pmatrix{
\tilde\ba_{0}-\bal_{0}(u)\vspace*{2pt}\cr
\tilde\ba-\bal(u)\vspace*{2pt}\cr
\tilde\bbe-\bbe_0
}
=-f^{-1}_U(u)\bD^{-1}_1(u)\bW^*_{n,1}(u)+o_p(1),
\end{equation}
where $\bW^*_{n,1}(u)=\frac{1}{\sqrt{nh}}\sumq\sumn K_i(u)\eta^*_{i,k}(u)(\be_k^T,\bX_i^T,\bZ_i^T)^T$.
Let
\[
\bW^\#_{n,1}(u)=\frac{1}{\sqrt{nh}}\sumq\sumn K_i(u)\eta_{i,k}(\be_k^T,\bX_i^T,\bZ_i^T)^T.
\]
Note that
\[
\operatorname{Cov}({\eta}_{i,k},{\eta}_{i,k'})={\tau}_{kk'},\qquad\operatorname{Cov}({\eta}_{i,k},\eta_{j,k'})=0\qquad\mbox{if }i\neq j.
\]
By some calculations, we have that $E[\bW^\#_{n,1}(u)]=\bzr$ and
$\operatorname{Var}[\bW^\#_{n,1}(u)]\rightarrow f_U(u)\nu_0\bSig_1(u)$. By the
Cram\'{e}r--Wold theorem,  the central limit theorem for $\bW_{n,1}(u)$
holds. Therefore,
%$(\bW_{n,1} -E[\bW_{n,1}])/\sqrt{Var[\bW_{n,1}]} \wcon N(\bzr, \bI).$
%It is easy to calculate that $Var[\bW_{n,1}] \pcon f_U(u) \bSig_1(u)$. Combined with
%(\ref{cltwn}), we have
\[
\bW^\#_{n,1}(u)\wcon N(\bzr,f_U(u)\nu_0\bSig_1(u)).
\]
Moreover, we have $ \operatorname{Var}[\bW^*_{n,1}(u)-\bW^\#_{n,1}(u)|U,\bX,\bZ]
%= \frac{1}{nh} \sumn K^2_i (\be_k^T,\bX_i^T,\bZ_i^T)^T (\be_k^T,\bX_i^T,\bZ_i^T) Var(\sumq \eta^*_{i,k} - \eta_{i,k})
\le\frac{q^2}{nh}\sumn K^2_i(u)(\be_k^T,\bX_i^T,\break\bZ_i^T)^T(\be_k^T,\bX_i^T,\bZ_i^T)\max_{k}\{F(c_k+|r_i|)-F(c_k)\}=o_p(1)$,
thus
\[
\operatorname{Var}[\bW^*_{n,1}(u)-\bW^\#_{n,1}(u)]=o(1).
\]
So, by Slutsky's theorem, conditioning on $\{U,\bX,\bZ\}$, we have
%e7.6 ###
\begin{equation}\label{CQRalphaR2}
\bW^*_{n,1}(u)-E[\bW^*_{n,1}(u)]\wcon N(\bzr,f_U(u)\nu_0\bSig_1(u)).
\end{equation}
We now calculate the conditional mean of $\bW^*_{n,1}(u)$:
\begin{eqnarray}\label{CQRalphaR3}
&&
\frac{1}{\sqrt{nh}} E[\bW^*_{n,1}(u) |U,\bX,\bZ]\nn
\\
&&\qquad=
\frac{1}{nh}\sumq\sumn K_i(u)\bigl\{F\bigl(c_k-r_{i}(u)\bigr)-F(c_k)\bigr\}(\be_k^t,\bX_i^T,\bZ_i^T)^T\nn
\\[-8pt]\\[-8pt]
&&\qquad=
-\frac{1}{nh}\sumq\sumn K_i(u)r_{i}(u)f(c_k)\{1+o(1)\}(\be_k^t,\bX_i^T,\bZ_i^T)^T\nn
\\
&&\qquad=
-\frac{\mu_2 h^2}{2}f_U(u)\bD_1(u)
\pmatrix{
 \bal''_{0}(u)\vspace*{3pt}\cr
 \bal''(u)\vspace*{2pt}\cr
 \bzr
}+o_p(h^2).\nn
\end{eqnarray}
The proof is completed by combining (\ref{CQRalphaR1}),
(\ref{CQRalphaR2}) and  (\ref{CQRalphaR3}).
\end{pf*}

\begin{pf*}{Proof of Theorem \ref{ThmCqrBeta1}}
%Let $v_k=\sqrt{n}(m_k-c_k)$ and $\bv=(v_1,\cdots,v_q)^T$.
Let $\bth=\sqrt{n}(\bbe-\bbe_0)$. Then,
\begin{eqnarray*}
&&
Y_i-\tilde a_{0k}(U_i)-\bX_i^T\tilde\ba(U_i)-\bZ_i^T\bbe
\\
%   &=& \alp_0(U_i) + \bX_i^T\bal(U_i) + \bZ_i^T\bbe_0 +\eps_i - m_{k} -  \tilde a_{0}(U_i) - \bX_i^T\tilde \ba(U_i)  - \bZ_i^T\bbe\\
&&\qquad=
\eps_i-c_k-\{\tilde a_{0k}(U_i)-\alp_0(U_i)-c_k\}-\bX_i^T\{\tilde\ba(U_i)-\bal(U_i)\}-\bZ_i^T(\bbe-\bbe_0)
\\
&&\qquad=
\eps_i-c_k-\tilde r_{i,k}-\bZ_i^T\bth/\sqrt{n},
%   &=& \eps_i - c_k  - \tilde r_i - (v_k+\bZ_{i}^T\bth)/\sqrt{n},
\end{eqnarray*}
where $\tilde r_{i,k}=\{\tilde a_{0k}(U_i)-\alp_0(U_i)-c_k\}+\bX_i^T\{\tilde\ba(U_i)-\bal(U_i)\}$. Then,
\[
\hat\bth=\argmin\sumq\sumn\rho_{\tau_k}\bigl(Y_i-\tilde a_{0k}(U_i)-\bX_i^T\tilde\ba(U_i)-\bZ_i^T\bbe\bigr)
\]
is also the minimizer of
\[
L_n(\bth)=\sumq\sumn\bigl\{\rho_{\tau_k}\bigl(\eps_i-c_k-\tilde r_{i,k}-\bZ_i^T\bth/\sqrt{n}\bigr)-\rho_{\tau_k}(\eps_i-c_k-\tilde r_{i,k})\bigr\}.
\]
By applying the identity (\ref{Knight}), we can rewrite $L_n(\bth)$ as
follows:
\begin{eqnarray*}
L_n(\bth)
&=&
\sumq\sumn\biggl\{\frac{\bZ_i^T\bth}{\sqrt{n}} [I(\eps_i \leq c_k)-\tau_k] \nn
\\
&&\hphantom{\sumq\sumn\biggl\{}{}+
\int_{\tilde r_{i,k}}^{\tilde r_{i,k}+\bZ_i^T\bth/\sqrt{n}} [  I(\eps_i \leq  c_k+z) - I(\eps_i \leq c_k) ]\,dz\biggr\}
\\
&=&
\Biggl(\frac{1}{\sqrt{n}}\sumq\sumn\eta_{i,k}\bZ_i\Biggr)^T\bth+B_{n}(\bth),
\end{eqnarray*}
where $B_{n}(\bth)=\sumq\sumn\int_{\tilde r_{i,k}}^{\tilde r_{i,k}+\bZ_i^T\bth/\sqrt{n}}[I(\eps_i\leq c_k+z)-I(\eps_i\leq c_k)]\,dz$.
Let us now calculate the conditional expectation of $B_{n}(\bth)$:
\begin{eqnarray*}
&&
E[B_{n}(\bth)|U,\bX,\bZ]
\\
%   &=& \sumq\sum_{i=1}^{n}
%         \int_{\tilde r_{i,k}}^{\tilde r_{i,k}+\bZ_i^T\bth/\sqrt{n}} [  F(c_k + z) - F(c_k)
%        ]   dz  \nonumber \\
&&\qquad=
\sumq\sum_{i=1}^{n}\int_{\tilde r_{i,k}}^{\tilde r_{i,k}+\bZ_i^T\bth/\sqrt{n}}[z f(c_k)\{1+o(1)\}]\,dz
\\
%   &=&  \sum_{i=1}^{n}
%         \Delta_{i,k}^2 f(c_k)/2  +o_p(1) \nonumber \\
&&\qquad=
\frac{1}{2}\bth^T\Biggl(\frac{1}{n}\sumq\sumn f(c_k)\bZ_i\bZ_i^T\Biggr)\bth-\Biggl(\frac{1}{\sqrt{n}}\sumq\sumn f(c_k)\tilde r_{i,k}\bZ_i\Biggr)^T\bth+o_p(1).
%   &&   + o_p(1) \\
\end{eqnarray*}
Define $R_{n}(\bth)=B_{n}(\bth)-E[B_{n}(\bth)|U,\bX,\bZ]$. It can
be shown that $R_{n}(\bth)=o_p(1)$. Hence,
%{\allowdisplaybreaks
\begin{eqnarray*}%\label{Ln4}
L_n(\bth)
&=&
\Biggl(\frac{1}{\sqrt{n}}\sumq\sumn\eta_{i,k}\bZ_i\Biggr)^T\bth+E[B_{n}(\bth)|U,\bX,\bZ]+R_{n}(\bth)
\\
&=&
\frac{1}{2}\bth^T\bS_n\bth+\Biggl(\frac{1}{\sqrt{n}}\sumq\sumn\eta_{i,k}\bZ_i\Biggr)^T\bth
-
\Biggl(\frac{1}{\sqrt{n}}\sumq\sumn f(c_k)\tilde r_{i,k}\bZ_i\Biggr)^T\bth
\\
&&{}+
o_p(1),
\end{eqnarray*}
where $\bS_n=\frac{1}{n}\sumq\sumn f(c_k)\bZ_i\bZ_i^T$. By
(\ref{CQRalphaR}), the third term in the previous expression can be
expressed as
\begin{eqnarray*}
%   &&  \frac{1}{\sqrt{n}} \sumq \sumn [ \eta_{i,k} -f(c_k)\tilde r_i ] \bZ_{i,k}^*   \\
&&
\frac{1}{\sqrt{n}}\sumq\sumn f(c_k)\tilde r_{i,k}\bZ_i
\\
&&\qquad=
\frac{1}{\sqrt{n}}\sumq\sumn\frac{f(c_k)}{f_U(U_i)}\bZ_i
\pmatrix{
\be_{k}\vspace*{2pt}\cr
\bX_{i}\vspace*{2pt}\cr
\bzr
}^T
\\
%   (\be_k^T ,\bX_{i}^T , \bzr)
&&\qquad\hphantom{=\frac{1}{\sqrt{n}}\sumq\sumn}{}\times
\bD_1^{-1}(U_{i})\left(\frac{1}{nh}\sumqq\sumnn\eta^*_{i',k'}(U_i)
\pmatrix{
\be_{k'}\vspace*{2pt}\cr
\bX_{i'}\vspace*{2pt}\cr
\bZ_{i'}
}
%(\be_{k'}^T,\bX_{i'}^T, \bZ_{i'}^T)^T
K_{i'}(U_i)\right)
\\
&&\qquad\quad{}+
O_p\bigl(h^{3/2}+\log^{1/2}(1/h)/\sqrt{nh^2}\bigr)
\\
&&\qquad=
\frac{1}{\sqrt{n}}\sumqq\sumnn\eta_{i',k'}\bolds{\delta}_{k'}(U_{i'},\bX_{i'},\bZ_{i'})+O_p\bigl(n^{1/2}h^{2}+\log^{1/2}(1/h)/\sqrt{nh^2}\bigr)
\\
&&\qquad=
\frac{1}{\sqrt{n}}\sumq\sumn\eta_{i,k}\bolds{\delta}_{k}(U_{i},\bX_{i},\bZ_{i})+o_p(1),
\end{eqnarray*}
where
%$\bolds{\xi}_{i,k}(U_{i}) = \sumqq \phi_{i,k,k'} f(c_{k'}) E[  \bZ | U=U_i ] + c \sum_{j=1}^{d_1} \phi_{i,k,q+j} E[ X_j \bZ | U=U_i ]$.
\[
\bolds{\delta}(U_{i},\bX_{i},\bZ_{i})=E[\bZ(\bc^T, c\bX^T,\mathbf{0})|U=U_i]\bD_1^{-1}(U_i)(I_q,\bone^T\bX_i,\bone^T\bZ_i)^T.
\]
%
%where $\bolds{\zeta}_{i,\tau}(U_i) = \xi_{i,1,\tau} E[ f_{\tau}(0|U,\bX,\bZ) \bZ | U=U_i ] + \sum_{j=1}^{d_1} \xi_{i,j+1,\tau} E[ f_{\tau}(0|U,\bX,\bZ) X_j \bZ | U=U_i ].$.
%
Therefore,
%{\allowdisplaybreaks
\begin{eqnarray*}%\label{Ln4}
L_n(\bth)
&=&
\frac{1}{2}\bth^T\bS_n\bth+\Biggl(\frac{1}{\sqrt{n}}\sumq\sumn\eta_{i,k}\{\bZ_i-\bolds{\delta}_{k}(U_{i},\bX_{i},\bZ_{i})\}\Biggr)^T\bth+o_p(1)
\\
&\triangleq&
\frac{1}{2}\bth^T\bS_n\bth+\bW_n^T\bth+o_p(1).
\end{eqnarray*}
It can be shown that $\bS_n=E(\bS_n)+o_p(1)=c\bS+o_p(1)$. Hence,
\[%\label{Ln5}
L_n(\bth)=\frac{c}{2}\bth^T\bS\bth+\bW_n^T\bth+o_p(1).
\]
Since the convex function $L_n(\bth)-\bW_n^T\bth$ converges in
probability to the convex function $\frac{c}{2}\bth^T\bS\bth$, it
follows, by the convexity lemma \cite{Pollard1991}, that the quadratic
approximation to $L_n(\bth)$ holds uniformly for $\bth$ in any compact
set $\Theta$. Thus, it follows that
%e7.7 ###
\begin{equation}
\hat{\bth}=-\frac{1}{c}\bS^{-1}\bW_n + o_p(1).
\end{equation}
By the Cram\'{e}r--Wold theorem, the central limit theorem for $\bW_n$
holds and  $\operatorname{Var}(\bW_n)\rightarrow\bDel=\sumq\sumqq\tau_{kk'}E\{\bZ-\bolds{\delta}_{k}(U,\bX,\bZ)\}
\{\bZ-\bolds{\delta}_{k'}(U,\bX,\bZ)\}^T$.
%
 %= \sumq \sumqq \tau_{kk'} E (\bZ_{1,k}^*-\boldsymbol{\gamma}_{1,k}) (\bZ_{1,k'}^*-\boldsymbol{\gamma}_{1,k'})^T$.
Therefore, the asymptotic normality of $\hat\bbe$ is followed by
%    \sqrt{n}( \hat\bth -\bth_0 ) \wcon
%    N(0, \frac{1}{c^2}\bS^{-1}\bXi\bS^{-1} ).
%Consequently,
\[
\sqrt{n}(\hat\bbe-\bbe_0)\wcon N\biggl(0,\frac{1}{c^2}\bS^{-1}\bDel\bS^{-1}\biggr).
\]
This completes the proof.
\end{pf*}

\begin{pf*}{Proof of Theorem \ref{ThmCqrAlpha}}
The asymptotic normality of $\hat\alp_0(u)$ and $\hat\bal(u)$ can be
obtained by following the ideas in the proof of Theorem \ref{ThmCqrAlpha1}.
\end{pf*}

\begin{pf*}{Proof of Theorem \ref{ThmCqrOracle}}
Use the same notation as in the proof of Theorem~\ref{ThmCqrBeta1}. Minimizing
\[
%    G_n(\bbe) =
\sumq\sumn\rho_{\tau_k}\{Y_i-\tilde a_{0k}(U_i)-\bX_i^T\tilde\ba(U_i)-\bZ_i^T\bbe\}+nq\sumd p'_{\lam_j}\bigl(\bigl|\beta_j^{(0)}\bigr|\bigr)|\beta_j|
\]
is equivalent to minimizing
\begin{eqnarray*}
G_n(\bth)
&=&
\sumq\sumn\bigl\{\rho_{\tau_k}\bigl(\eps_i-c_k-\tilde r_{i,k}-\bZ_{i}^T\bth/\sqrt{n}\bigr)-\rho_{\tau_k}(\eps_i-c_k-\tilde r_{i,k})\bigr\}
\\
&&{}+
nq\sumd p'_{\lam_j}\bigl(\bigl|\beta_j^{(0)}\bigr|\bigr)(|\beta_j|-|\beta_{0j}|)
\\
%    &=&L_n(\bth)+n \sumd p'_{\lam_j}(|\beta_j^{(0)}|)(|\beta_j|-|\beta_{0j}|)\\
&=&
\frac{c}{2}\bth^T\bS\bth+\bW_n^T\bth+nq\sumd p'_{\lam_j}\bigl(\bigl|\beta_j^{(0)}\bigr|\bigr)(|\beta_j|-|\beta_{0j}|)+ o_p(1),
\end{eqnarray*}
where $\bth=\sqrt{n}(\bbe-\bbe_0)$ and $\tilde r_{i,k}=\{\tilde a_{0k}(U_i)-\alp_0(U_i)-c_k\}+\bX_i^T\{\tilde\ba(U_i)-\bal(U_i)\}$.
Similar to the derivation in the proof of Theorem 5 in
Zou and Li \cite{ZouLi2008}, the third term above can be expressed as
%e7.8 ###
\begin{equation}\label{results3}
nq\sumd p'_{\lam_j}\bigl(\bigl|\beta_j^{(0)}\bigr|\bigr)(|\beta_j|-|\beta_{0j}|)
\pcon
\cases{
0,&\quad if $\bbe_{2}=\bbe_{20}$,\cr
\infty,&\quad otherwise.}
\end{equation}
Therefore, by the epiconvergence results \cite{Geyer1994,KnightFu2000}, we have $\hat \bbe{}^{\mathrm{OSE}}_2 \pcon 0$
and the asymptotic results for $\hat\bbe{}^{\mathrm{OSE}}_1$ holds.

To prove sparsity, we only need to show that $\hat\bbe{}^{\mathrm{OSE}}_2=0$
with probability tending to 1. It  suffices to prove that if
$\beta_{0j}=0$, then  $P(\hat\beta_j^{\mathrm{OSE}}\neq 0 )\rar 0$. By using
the fact that $|\frac{\rho_{\tau}(t_1)-\rho_{\tau}(t_2)}{t_1-t_2}|\leq\max(\tau,1-\tau)<1$, if $\hat\beta_j^{\mathrm{OSE}}\neq 0$, then we must have
$\sqrt{n} p'_{\lam_j}(|\beta_j^{(0)}|)<\frac{1}{n}\sumn |Z_{ij}|$.
Thus, we have $P(\hat\beta_j^{\mathrm{OSE}}\neq 0)\leq P(\sqrt{n}p'_{\lam_j}(|\beta_j^{(0)}|)<\frac{1}{n}\sumn|Z_{ij}|)$. However,
under the assumptions, we have $\sqrt{n}p'_{\lam_j}(|\beta_j^{(0)}|)\rar\infty$. Therefore, $P(\hat\beta_j^{\mathrm{OSE}}\neq 0)\rar 0$. This
completes the proof.
%\rightqed
\end{pf*}

\printaddresses

\end{document}